\documentclass[12pt]{article}

\usepackage{amsfonts}
\usepackage{amsmath}
\usepackage{enumerate}
\usepackage{tabu}
\usepackage{yfonts}

\newcommand\F{\mathbb{F}}
\newcommand\Z{\mathbb{Z}}

\newcommand\cF{\mathcal{F}}
\newcommand\cH{\mathcal{H}}
\newcommand\cJ{\mathcal{J}}

\newcommand\al{\alpha}
\newcommand\bt{\beta}
\newcommand\gm{\gamma}
\newcommand\dl{\delta}

\newcommand\lm{\lambda}
\newcommand\sg{\sigma}
\newcommand\om{\omega}

\newcommand\Lm{\Lambda}

\newcommand\fa{\textfrak{a}}
\newcommand\fb{\textfrak{b}}
\newcommand\fc{\textfrak{c}}
\newcommand\fd{\textfrak{d}}
\newcommand\fe{\textfrak{e}}
\newcommand\fs{\textfrak{s}}
\newcommand\ft{\textfrak{t}}

\newcommand\fC{\textfrak{C}}
\newcommand\fD{\textfrak{D}}
\newcommand\fJ{\textfrak{J}}
\newcommand\fS{\textfrak{S}}
\newcommand\fT{\textfrak{T}}
\newcommand\fW{\textfrak{W}}

\newcommand\bsg{\boldsymbol{\sg}}

\newcommand\la{\langle}
\newcommand\ra{\rangle}
\newcommand\lla{\la\!\la}
\newcommand\rra{\ra\!\ra}
\newcommand\ad{\mathrm{ad}}
\newcommand\im{\mathrm{im}\,}

\newcommand\Aut{\mathrm{Aut}}
\newcommand\Miy{\mathrm{Miy}}

\newtheorem{main}{Theorem}
\newtheorem{thm}{Theorem}[section]
\newtheorem{lem}[thm]{Lemma}
\newtheorem{prop}[thm]{Proposition}
\newtheorem{cor}[thm]{Corollary}
\newtheorem{defn}[thm]{Definition}
\newtheorem{exa}[thm]{Example}
\newtheorem{que}[thm]{Question}
\newtheorem{conj}{Conjecture}

\newcommand\pf{\noindent\textbf{Proof:\quad}}
\newcommand\qed{\hfill$\Box$}

\title{Solid subalgebras in algebras of Jordan type half}
\author{I.~Gorshkov, S.~Shpectorov, and A.~Staroletov}

\date{\vspace{-30px}}

\begin{document}
	\maketitle
	
	\begin{abstract}
		The class of algebras of Jordan type $\eta$ was introduced by Hall, Rehren 
		and Shpectorov in 2015 within the much broader class of axial algebras. 
		Algebras of Jordan type are commutative algebras $A$ over a field of 
		characteristic not $2$, generated by primitive idempotents, called axes, whose 
		adjoint action on $A$ has minimal polynomial dividing $(x-1)x(x-\eta)$ 
		and where multiplication of eigenvectors follows the rules similar to 
		the Peirce decomposition in Jordan algebras. 
		
		Naturally, Jordan algebras generated by primitive idempotents are examples 
		of algebras of Jordan type $\eta=\frac{1}{2}$. Further examples are given 
		by the Matsuo algebras constructed from $3$-transposition groups. These examples 
		exist for all values of $\eta\neq 0,1$. Jordan algebras and (factors of) Matsuo 
		algebras constitute all currently known examples of algebras of Jordan type and 
		it is conjectured that there are now additional examples.
		
		In this paper we introduce the concept of a solid $2$-generated subalgebra, as 
		a subalgebra $J$ such that all primitive idempotents from $J$ are axes of $A$. 
		We prove that, for axes $a,b\in A$, if $(a,b)\notin\{0,\frac{1}{4},1\}$ then 
		$J=\lla a,b\rra$ is solid, that is, generic $2$-generated subalgebras are solid. 
		Furthermore, in characteristic zero, $J$ is solid even for the values $(a,b)=0,1$.
		As a corollary, in characteristic zero, either $A$ has infinitely many axes and 
		an infinite automorphism group, or it is a Matsuo algebra or a factor of Matsuo 
		algebra. 
	\end{abstract}
	
	
	{\bf MSC classes:} 17A99; 20F29

	\newcommand{\Addresses}{{
			\bigskip\noindent
			\footnotesize
			Ilya~Gorshkov, \textsc{Sobolev Institute of Mathematics, Novosibirsk, Russia;}\\\nopagebreak
			\textit{E-mail address: } \texttt{ilygor8@gmail.com}
			
			\medskip\noindent
			Sergey~Shpectorov, \textsc{School of Mathematics, University of Birmingham,
				Watson Building, Edgbaston, Birmingham B15 2TT, United Kingdom;}
			\\\nopagebreak
			\textit{E-mail address: } \texttt{s.shpectorov@bham.ac.uk}
			
			\medskip\noindent
			Alexey~Staroletov, \textsc{Sobolev Institute of Mathematics, Novosibirsk, Russia;}\\\nopagebreak
			\textit{E-mail address: } \texttt{staroletov@math.nsc.ru}
			
			\medskip
	}}
	
	\section{Introduction} \label{introduction}
	
	The class of algebras of Jordan type $\eta$ was introduced by Hall, Rehren 
	and Shpectorov in \cite{hrs2} within the much broader class of axial algebras\footnote{Axial algebras were introduced by Hall, Rehren and Shpectorov in \cite{hrs1} and \cite{hrs2} broadly generalising some features of 
		Majorana algebras of Ivanov \cite{i}, which in turn were modelled after the key properties of the Griess algebra \cite{g,c} for the Monster sporadic simple group. The area of axial algebras is currently experiencing rapid growth with new connections found within many areas of mathematics and beyond (cf. \cite{ms}).}. 
	As such, an algebra of Jordan type $\eta$ is a commutative non-associative algebra 
	$A$ over a field $\F$, such that $A$ is generated by a set of special idempotents 
	called (primitive) axes. The special property of axes is that they obey a fusion 
	law. For the class of algebras of Jordan type, the relevant fusion law is shown in 
	Figure \ref{Jordan law}.
	\begin{figure}[ht] \label{Jordan law}
		\begin{center}
			\renewcommand{\arraystretch}{1.3}
			\begin{tabular}{|c||c|c|c|}
				\hline
				$\ast$ & $1$ & $0$ & $\eta$\\
				\hline\hline
				$1$ & $1$ & & $\eta$\\
				\hline
				$0$ & & $0$ & $\eta$\\
				\hline
				$\eta$ & $\eta$ & $\eta$ & $1,0$\\
				\hline
			\end{tabular}
			\caption{Jordan type fusion law $\cJ(\eta)$}
		\end{center}
	\end{figure}
	
	For the exact meaning of this and other fusion laws and the definition of primitivity 
	see Section \ref{background}. However, at the ideological level, the above condition 
	means that, for every axis, the algebra admits a decomposition similar to Peirce 
	decompositions in Jordan algebras, and this explains the name given to this class of 
	algebras.
	
	Throughout, we will assume that $\F$ is a field of characteristic not $2$, although 
	the definition itself makes sense even in this case.
	
	Jordan algebras generated by primitive idempotents are, of course, examples of algebras 
	of Jordan type, but they only arise for $\eta=\frac{1}{2}$. Further examples are given 
	by the Matsuo algebras introduced by Matsuo in \cite{m} (see also a broader definition 
	in \cite{hrs2}), which are constructed from $3$-transposition groups and exist for all 
	values of $\eta\neq 1,0$. Jordan algebras and (factors of) Matsuo algebras constitute all 
	currently known examples of algebras of Jordan type and in fact it is conjectured (e.g. see \cite{ms}) that no further examples exist:
	
	\begin{conj}
		Every algebra of Jordan type is either a Jordan algebra, or a Matsuo algebra, or a factor of a Matsuo algebra.\footnote{Unlike Jordan algebras, a factor of a Matsuo algebra does not need to 
			be a Matsuo algebra, but it remains an algebra of Jordan type. Hence we should include, 
			together with Matsuo algebras, all their non-zero factors.}    
	\end{conj}
	
	This conjecture can be considered as a generalisation to the axial context of the classification of Jordan algebras. While this conjecture remains open, some of its cases have been completed. In the same paper \cite{hrs2} (with a minor correction by Hall, 
	Segev and Shpectorov \cite{hss1}), it was shown that, when $\eta\neq\frac{1}{2}$, every 
	algebra of Jordan type $\eta$ is isomorphic to a Matsuo algebra or a factor of a Matsuo 
	algebra. Thus, only the case of 
	$\eta=\frac{1}{2}$ remains open. The difficulty of this case stems from the fact that 
	we have here two very different classes of examples, Matsuo algebras for 
	$\eta=\frac{1}{2}$ and Jordan algebras. 
	
	Partial results were also obtained for $\eta=\frac{1}{2}$. In the original paper 
	\cite{hrs2}, all $2$-generated algebras of Jordan type were classified. In 
	\cite{hss2}, Hall, Segev, and Shpectorov showed that every algebra of Jordan type is 
	metrisable, that is, it admits a non-zero bilinear form associating with the algebra 
	product (so-called Frobenius form). Additionally, this form is unique if we require that 
	$(a,a)=1$ for every generating axis of $A$. Gorshkov and Staroletov \cite{gs} completely 
	classified $3$-generated algebras of Jordan type. In particular, for $\eta=\frac{1}{2}$, 
	they showed that there exists a $4$-parameter family of $3$-generated algebras of Jordan 
	type half, which are all of dimension $9$. All other $3$-generated algebras of Jordan type 
	half are factors of these $9$-dimensional algebras. Furthermore, they proved that all 
	$3$-generated examples are Jordan algebras, which means that the families of Matsuo algebras and 
	Jordan algebras split from each other when the number of generators is at least $4$. Note 
	that the intersection of these two classes was determined by De Medts and Rehren 
	\cite{dmr}. Their result was recently enhanced by Gorshkov, Mamontov, and Staroletov 
	\cite{gms}, who also classified, only in characteristic zero, all factors of Matsuo 
	algebras that are Jordan algebras. Finally, De Medts, Rowen, and Segev \cite{dmrs} 
	recently proved that every $4$-generated algebra of Jordan type half has dimension at most 
	$81$. We are still far from the classification of $4$-generated algebras of Jordan type, but it is very reassuring that the dimension in this case
	is bounded.
	
	In this paper we look into the structure of an arbitrary algebra of Jordan type half. As 
	we already mentioned, all algebras of Jordan type half generated by two axes, say $a$ and 
	$b$, are known. With a finite number of exceptions, the algebra $\lla a,b\rra$\footnote{In 
		an algebra, we use angular brackets $\la~\ra$ for the linear span of a set, while double 
		brackets $\lla~\rra$ stand for subalgebra generation.} is $3$-dimensional determined by a 
	single parameter $\al=(a,b)$, the value of the Frobenius form on the pair $a$ and $b$. 
	Furthermore, for almost all values of $\al$, the algebra $\lla a,b\rra$ is the same up to 
	isomorphism if we additionally assume that the field is algebraically closed. Hence the 
	issue is not what the algebra is, but rather how the pair $a$, $b$ embeds into it, and 
	this is determined by the value of $\al$.
	
	In a general algebra $A$ of Jordan type half, we can consider the family of all subalgebras
	generated by two axes. Generically, a $2$-generated subalgebra $\lla a,b\rra\subseteq A$ 
	contains many more (infinitely many, if the field is infinite) primitive axes. A natural 
	question is: can we somehow check whether these additional primitive idempotents are axes 
	in the entire $A$? We say that $\lla a,b\rra$ is \emph{solid} if, in fact, all its 
	primitive idempotents are axes in $A$. The surprising result that we prove in this paper 
	is that almost all values of $\al=(a,b)$ 
	guarantee that $\lla a,b\rra$ is solid. 
	
	\begin{main} \label{main solid}
		Suppose that $A$ is an algebra of Jordan type half and $J=\lla a,b\rra\subseteq A$, where 
		$a$ and $b$ are primitive axes from $A$. If $\al=(a,b)\notin\{0,\frac{1}{4},1\}$ then $J$ 
		is solid, that is, all idempotents $c\in J$, $c\neq 0,1_J$, are primitive axes in $A$.
	\end{main}
	
	This gives us a much richer set of primitive axes in $A$ than we could ever expect. 
	Furthermore, over the fields of characteristic zero, additional 
	arguments borrowed from algebraic geometry give us even more.
	
	\begin{main} \label{exceptional}
		Let $A$ be an algebra of Jordan type half over a filed $\F$ of characteristic zero. If 
		$a,b\in A$ are primitive axes with $(a,b)=0$ or $1$ then $J=\lla a,b\rra$ is solid.
	\end{main}
	
	Hence, in characteristic zero all non-solid $2$-generated subalgebras $J=\lla a,b\rra$ 
	arise only in the case where $(a,b)=\frac{1}{4}$. In this case, $J\cong 3C(\frac{1}{2})$
	and it contains exactly three primitive axes of $A$, namely $a$, $b$, and $c=a^{\tau_b}=
	b^{\tau_a}$, where $\tau_a$ and $\tau_b$ are the Miyamoto involutions (automorphisms of 
	$A$ of order $2$) associated with $a$ and $b$. (See Section \ref{background} for the 
	exact definitions.)
	
	We note that over infinite fields $\F$ all solid lines, except the ones where $ab=0$
	and so $J=\lla a,b\rra\cong \F\oplus\F=2B$, contain infinitely many primitive axes. 
	The algebra $2B$ contains only two primitive axes, $a$ and $b$.
	This means that if the algebra $A$ contains only finitely many axes (equivalently, 
	has a finite automorphism group) then all $2$-generated subalgebras in $A$ are of 
	type $2B$ or $3C$ and this immediately leads to the following striking result.
	
	\begin{main} \label{main corollary}
		If $\F$ is a field of characteristic zero and $A$ is an algebra of Jordan type half over 
		$\F$ then either $A$ contains infinitely many primitive axes or $A$ is a Matsuo algebra 
		or a factor of Matsuo algebra.
	\end{main}
	
	This provides a delineation between the infinite case where we expect the Jordan algebras 
	to be the principal examples and the finite case, where we can only have the Matsuo 
	algebras and their factors. Note that this is not a strict separation of cases: not only 
	some smaller Matsuo algebras are Jordan, but also, as we show in Example \ref{vertical} 
	in this paper, 
	Matsuo algebras can be non-Jordan and still contain solid $2$-generated subalgebras with 
	infinitely many primitive axes. Still, we believe that we found an important structural 
	property of algebras of Jordan type half, which suggests that the focus should be on the 
	geometry of all primitive axes. 
	
	Furthermore, note that the non-Jordan algebra in Example \ref{vertical} have only part 
	of its $2$-generated subalgebras solid. This allows us to formulate the following.
	
	\begin{conj}
		If $A$ is an algebra of Jordan type half over a field $\F$ with $\mathrm{char}\,\F=0$, 
		such that all $2$-generated subalgebras in it are solid, then $A$ is a Jordan algebra.
	\end{conj}
	
	Let us now describe the contents of this paper. Following this introduction, in Section 
	\ref{background}, we provide the necessary background on axes and axial algebras. In 
	Section \ref{2-gen}, we discuss $2$-generated algebras of Jordan type half, as classified 
	in \cite{hrs2}, and establish some properties that are used later in the paper. Similarly, 
	in Section \ref{3-gen}, based on \cite{gs}, we provide a description of all $3$-generated 
	algebras of Jordan type half and establish some further properties. As the description of 
	algebras in this section is quite complicated, some of the statements are verified using 
	the computer algebra system GAP \cite{GAP}. In Section \ref{solid}, we introduce solid 
	subalgebras and prove our principal result, Theorem \ref{main solid}, stating that almost 
	all $2$-generated subalgebras are solid. In Section \ref{char zero}, we restrict ourselves 
	to the case of the field of characteristic zero and prove that two of the ``exceptional'' 
	values of the parameter $\al=(a,b)$ still lead, in this case, to solid subalgebras, i.e., 
	they are non-exceptional. Finally, in Section \ref{example}, we show an example of a 
	non-Jordan Matsuo algebra, which admits solid subalgebras with infinitely many primitive 
	axes. We also include there a number of open questions. 
	
	\medskip
	This work was supported by the Mathematical Center in Akademgorodok under agreement No. 
	075-15-2022-281 with the Ministry of Science and Higher Education of the Russian Federation. 
	
	\section{Background} \label{background}
	
	We start by discussing fusion laws.
	
	\begin{defn}
		A \emph{fusion law} is a pair $(\cF,\ast)$ where $\cF$ is a set and $\ast:\cF\times\cF\to 
		2^\cF$ is a binary operation on $\cF$ with values in the set $2^\cF$ of all subsets of 
		$\cF$. 
	\end{defn} 
	
	An example of a fusion law, the main one that we consider in this paper, is shown in 
	Figure \ref{Jordan law}. Here $\cF=\cJ(\eta)=\{1,0,\eta\}$ is a subset of a field $\F$ 
	and the binary operation $\ast$ is represented by a table similar to the group 
	multiplication table. Every cell of the table represents a specific set $\mu\ast\nu$ for 
	$\mu,\nu\in\cF$. Note that we simply list all elements of $\mu\ast\nu$ in the cell, 
	dropping the set brackets $\{\}$. In particular, the empty cell in the table represents 
	the empty set $\mu\ast\nu$.
	
	To simplify the notation, we will speak of the set $\cF$ as the fusion law, with the 
	binary operation being assumed.
	
	We use fusion laws to axiomatise classes of algebras. Here by an algebra we mean an 
	arbitrary non-associative algebra over a field $\F$, that is, we do not assume 
	associativity. 
	
	Suppose that $A$ is a commutative algebra over $\F$. Let us fix a fusion law 
	$\cF\subseteq\F$. For $a\in A$, we denote by $\ad_a$ the adjoint map $A\to A$ sending 
	$u$ to $au$. For $\lm\in\F$, we use $A_\lm(a)=\{u\in A\mid au=\lm u\}$ to denote the 
	$\lm$-eigenspace of $\ad_a$. Note that $A_\lm(a)\neq 0$ if and only if $\lm$ is an 
	eigenvalue of $\ad_a$. Furthermore, we extend the eigenspace notation to sets of numbers: 
	if $\Lm\subseteq\F$ then $A_\Lm(a)=\oplus_{\lm\in\Lm}A_\lm(a)$. 
	
	\begin{defn}
		We say that $0\neq a\in A$ is an ($\cF$)-axis if the following hold:
		\begin{enumerate}[$(1)$]
			\item $a$ is an idempotent, i.e., $a^2=a$;
			\item $\ad_a$ is semi-simple with all its eigenvalues in $\cF$; that is, $A=A_\cF(a)$;
			\item $a$ obeys the fusion law $\cF$: for $\lm,\mu\in\cF$, we have that 
			$A_\lm(a)A_\mu(a)\subseteq A_{\lm\ast\mu}(a)$.
		\end{enumerate}
	\end{defn}
	
	Hence the fusion law specifies the possible eigenvalues of the adjoint map of an axis and 
	limits how eigenvectors of $\ad_a$ can multiply.
	
	We note that we must have $1\in\cF$ since $a$ is an idempotent.
	
	\begin{defn}
		An axis $a\in A$ is called \emph{primitive} if $A_1(a)=\la a\ra$ is $1$-dimensional.
	\end{defn}
	
	Finally, we can introduce the class of axial algebras.
	
	\begin{defn}
		An \emph{($\cF$-)axial algebra} is a pair $(A,X)$ where $A$ is a commutative 
		non-associative algebra over $\F$ and $X\subset A$ is a set of ($\cF$-)axes generating 
		$A$, i.e., $A=\lla X\rra$. We say that $(A,X)$ is a \emph{primitive} axial algebra when 
		all axes in $X$ are primitive.
	\end{defn}
	
	We note that $X$ is just a generating set of axes and it is never assumed that $X$ 
	consists of all idempotents from $A$ having the properties of an axis. Depending on the 
	situation, it might be beneficial sometimes to change the generating set $X$ for another 
	set. Accordingly, we will say that $A$ is an axial algebra, meaning that there is some 
	generating set $X$ for it.
	
	In this paper we focus on the class of algebras of Jordan type $\eta=\frac{1}{2}$, which 
	by definition are primitive axial algebras for the fusion law $\cJ(\frac{1}{2})$ from 
	Figure \ref{Jordan type half}.
	\begin{figure}[ht] \label{Jordan type half}
		\begin{center}
			\renewcommand{\arraystretch}{1.5}
			\begin{tabular}{|c||c|c|c|}
				\hline
				$\ast$ & $1$ & $0$ & $\frac{1}{2}$\\
				\hline\hline
				$1$ & $1$ & & $\frac{1}{2}$\\
				\hline
				$0$ & & $0$ & $\frac{1}{2}$\\
				\hline
				$\frac{1}{2}$ & $\frac{1}{2}$ & $\frac{1}{2}$ & $1,0$\\
				\hline
			\end{tabular}
			\caption{Jordan type half}
		\end{center}
	\end{figure}
	
	This is the last open case of the more general class of algebras of Jordan type $\eta$, 
	$\eta\neq 1,0$, which is specified by the fusion law $\cJ(\eta)$ from Figure 
	\ref{Jordan law}.
	
	Known examples of algebras of Jordan type include Jordan algebras, which are commutative algebras $A$
	satisfying the \emph{Jordan identity}
	$$u(vu^2)=(uv)u^2,$$
	for all $u,v\in A$. While all idempotents in a Jordan algebra satisfy the fusion law $\cJ(\frac{1}{2}$
	(this follows from the Peirce decomposition result), not every Jordan algebra contains non-zero idempotents, 
	and so only the Jordan algebras generated by primitive idempotents are algebras of Jordan type half.
	
	The other know set of examples of algebras of Jordan type are Matsuo algebras.
	
	\begin{defn}
		Suppose that $(G,D)$ is a $3$-transposition group, i.e., $D$ is a normal set of involutions (elements of 
		order $2$ from $G$) such that $G=\la D\ra$ and $|cd|\leq 3$ for all $c,d\in D$. Let $\F$ be a field of 
		characteristic not $2$ and $\eta\in\F$, $\eta\neq 0,1$. The \emph{Matsuo algebra} $M=M_\eta(G,D)$ over 
		$\F$ has $D$ as its basis and the algebra product is defined on $D$ as follows:
		$$c\cdot d=\left\{
		\begin{array}{rl}
			c,&\mbox{if $c=d$;}\\
			0,&\mbox{if $|cd|=2$;}\\
			\frac{\eta}{2}(c+d-e),&\mbox{if $|cd|=3$,}
		\end{array}
		\right.$$
		where in the last case $c=a^b=b^a$.
	\end{defn}
	Related to the $3$-transposition group $(G,D)$ is its \emph{Fischer space}, which is a point-line 
	geometry, whose point set is $D$ and whose lines are all triples $\{c,d,e\}$, where $|cd|=3$ and 
	$e=c^d$. Clearly, the Matsuo algebra $M$ can also be defined in terms of the Fischer space. 
	
	In Matsuo algebras, the basis elements from $D$ are primitive axes of Jordan type $\eta$ and 
	so all Matsuo algebras and their factors are examples of algebras of Jordan type.
	
	\medskip
	Let us now mention some known properties of algebras of Jordan type.
	
	\begin{defn}
		A fusion law $\cF$ is called \emph{Seress} if $0\in\cF$ and 
		$1\ast\lm,0\ast\lm\subseteq\{\lm\}$ for all $\lm\in\cF$. Axial algebras with a Seress
		fusion law are also said to be Seress.
	\end{defn}
	
	Note that $\cJ(\eta)$ is Seress, and so algebras of Jordan type are Seress.
	
	\begin{prop}[Seress Lemma] \label{Seress Lemma}
		Suppose that $a\in A$ is an axis for a Seress fusion law then $a$ associates with
		$A_1(a)\oplus A_0(a)$, that is,
		$$a(vw)=(av)w$$
		for all $v\in A$ and $w\in A_1(a)\oplus A_0(a)$.
	\end{prop}
	
	We note that here $a$ is an arbitrary axis, not necessarily primitive.
	
	\begin{defn}
		A \emph{Frobenius form} on an axial algebra $A$ is a non-zero bilinear form that 
		associates with  the algebra product:
		$$(uv,w)=(u,vw)$$
		for all $u,v,w\in A$.
	\end{defn}
	
	Note that the Frobenius Form on an axial algebra is automatically symmetric.
	
	According to \cite{hss2}, every algebra $A$ of Jordan type admits a Frobenius form, that 
	is, it is \emph{metrisable}. Furthermore, this form is unique subject to the condition 
	that $(a,a)=1$ for each axis $a\in X$. (In fact, the same is true for all primitive axes 
	in $A$, regardless of whether or not they are included in the generating set $X$.) In the 
	remainder of the paper, we will always assume that the Frobenius form on $A$ satisfies 
	this additional condition.
	
	Note that, for an axis $a$, the eigenspaces $A_\lm(a)$ for different $\lm$ are pairwise 
	orthogonal with repect to the Frobenius form. In particular, it is easy to see that, when 
	$a$ is primitive, the orthogonal projection of an arbitrary $u\in A$ to the eigenspace 
	$A_1(a)$ is equal to $(u,a)a$.
	
	One can define a graph on the set $X$, called the \emph{projection graph}, where axes 
	$a$ and $b$ are adjacent whenever $(a,b)\neq 0$. We say that $A$ is \emph{connected} when 
	its projection graph is connected.
	
	Note that the radical $A^\perp$ of the Frobenius form on $A$ is an ideal, and it was 
	shown in \cite{kms} that, in a connected algebra $A$ of Jordan type, the radical 
	$A^\perp$ is the unique maximal proper ideal. In particular, a connected $A$ is simple if 
	and only if $A^\perp=0$, i.e., the Frobenius form on $A$ is non-degenerate.
	
	The fusion law of Jordan type is $C_2$-graded. (See \cite{dmpsvc} for a more general, 
	categorical point of view.) Indeed, if, for an axis $a\in A$, we define $A_+:=A_1(a)
	\oplus A_0(a)$ and $A_-:=A_\eta(a)$, then we obtain a $C_2$-grading of the algebra, 
	$A=A_+\oplus A_-$, where $A_+A_+,A_-A_-\subseteq A_+$ and $A_+A_-\subseteq A_-$. 
	Consequently, the linear map $\tau_a:A\to A$ that is equal to the identity map on $A_+$ 
	and minus the identity map on $A_-$ is an automorphism of $A$, called the \emph{Miyamoto 
		involution} (or, sometimes, tau involution). Indeed, $\tau_a^2=1$, although it is also 
	possible that $\tau_a=1$, which happens exactly when $A_-=0$. (Then $a$ obeys a more 
	simple fusion law and it can be shown that it associates with the entire $A$.)
	
	\begin{defn}
		The Miyamoto group $\Miy(A)$ of $A$ is the subgroup of the full automorphism group of 
		$A$ generated by the Miyamoto involutions $\tau_a$, $a\in X$.
	\end{defn}
	
	We note that if $a$ is an axis and $\varphi$ an automorphism then $a^\varphi$ is also 
	an axis. In particular, we can \emph{close} the set of generators $X$ by taking as the 
	new set of generators $\bar X=\{a^\varphi\mid a\in X,\varphi\in\Miy(A)\}$. It was shown 
	in \cite{kms} that the radical of $A$ and its Miyamoto group do not change when we close 
	the set of generating axes. Equally, other important properties of $A$ are unaffected 
	by this operation. 
	
	Later on we will need the following result from \cite{hss1}.
	
	\begin{thm} \label{basis of axes}
		If $(A,X)$ is an algebra of Jordan type and $X=\bar X$ is closed then $X$ spans $A$, that 
		is, $A$ admits a basis consisting of axes.
	\end{thm}
	
	Before we prove the main results of this article, we need some information on $2$- and 
	$3$-generated algebras of Jordan type half.

	\section{Two axes} \label{2-gen}
	
	This section is written largely based on \cite{hrs2}, except we use a 
	different notation for the algebras.
	
	Let $\fJ=\fJ(\al)$, $\al\in\F$, be the algebra over $\F$ with basis $\fa$, $\fb$, $\bsg$, and 
	multiplication given by the following table:
	$$
	\renewcommand{\arraystretch}{1.7}
	\begin{array}{|c|ccc|}
		\hline
		&\fa&\fb&\bsg\\
		\hline
		\fa&\fa&\frac{1}{2}\fa+\frac{1}{2}\fb+\bsg&\frac{\al-1}{2}\fa\\
		\fb&\frac{1}{2}\fa+\frac{1}{2}\fb+\bsg&\fb&\frac{\al-1}{2}\fb\\
		\bsg&\frac{\al-1}{2}\fa&\frac{\al-1}{2}\fb&\frac{\al-1}{2}\bsg\\
		\hline
	\end{array}
	$$
	
	This is an algebra of Jordan type half for all values of $\al$, with 
	$\fa$ and $\fb$ being primitive axes, and in fact, $\fJ$ is universal, as 
	the following result shows.
	
	\begin{thm}
		Let $A=\lla a,b\rra$ be an algebra of Jordan type half, where $a$ and $b$ 
		are primitive axes. Let $\al=(a,b)$. Then there is a unique surjective 
		homomorphism from $\fJ(\al)$ onto $A$ sending $\fa$ and $\fb$ onto $a$ and $b$, 
		respectively.
	\end{thm} 
	
	In particular, this result clarifies the meaning of the parameter $\al$: it 
	is equal to the value of Frobenius form: $\al=(a,b)=(\fa,\fb)$. The entire Gram matrix 
	of the Frobenius form on $\fJ$ with respect to the basis $\fa,\fb,\bsg$ 
	is as follows:
	$$
	\renewcommand{\arraystretch}{1.7}
	\left(
	\begin{array}{ccc}
		1&\al&\frac{\al-1}{2}\\
		\al&1&\frac{\al-1}{2}\\
		\frac{\al-1}{2}&\frac{\al-1}{2}&\frac{(\al-1)^2}{2}
	\end{array}
	\right),
	$$
	and its determinant is $-\frac{\al(\al-1)^3}{2}$. In particular, $\fJ$ is simple unless $\al=0$ or $1$. 
	
	\begin{thm} \label{all ideals}
		The algebra $\fJ(\al)$ is simple if and only if $\al\notin\{0,1\}$.
		Furthermore, 
		\begin{enumerate}[\rm(a)]
			\item $\fJ(0)$ has a unique minimal non-zero ideal 
			$I=\la\frac{1}{2}\fa+\frac{1}{2}\fb+\bsg\ra$ and the quotient 
			$\overline{\fJ(0)}=\fJ(0)/I$ is the algebra $\lla\bar\fa\rra\oplus\lla\bar 
			\fb\rra\cong\F\oplus\F$ (i.e., $\bar\fa\bar\fb=0$);
			\item $\fJ(1)$ is baric having exactly two proper non-zero ideals: 
			\begin{enumerate}[\rm(i)]
				\item $R=\la\fa-\fb,\bsg\ra$, the baric radical; and
				\item $Z=\la\bsg\ra$, the algebra annihilator; that is, $Z\fJ(1)=0$. 
			\end{enumerate}
		\end{enumerate}
	\end{thm}
	
	Recall that an algebra $A$ over $\F$ is called \emph{baric} if it admits 
	a surjective homomorphism $wt:A\to\F$, called the \emph{baric weight}. The 
	kernel of $wt$ is known as the \emph{baric radical} of $A$. Note that, for 
	every idempotents $a\in A$, we have that $wt(a)=0$ or $1$, and in the axial 
	case we additionally require that $wt(a)=1$ for all generating axes $a\in A$.
	That is, we do not want axes to be in the baric radical. In this sense, 
	we do not view the algebra $\overline{\fJ(0)}$ from Theorem \ref{all 
		ideals}(b)(i) (also known as $2B$) as a baric algebra, even though it has 
	two homomorphisms onto $\F$: either of the two sends one of the axes, 
	$\bar\fa$ or $\bar\fb$, to $0$.
	
	Note also, that $(u,v)=wt(u)wt(v)$ is a Frobenius form on a baric algebra, 
	and so the additional assumption above amounts to asking that all axes $a$ 
	in a baric axial algebra satisfy $(a,a)=1$. Also, the radical of the algebra
	coincides then with the baric radical.
	
	The final comment on Theorem \ref{all ideals} is that the $2$-dimensional 
	factor algebra $\overline{\fJ(1)}=\fJ(1)/Z$ is spanned by $\bar\fa$ and $\bar\fb$, 
	with the product given by $\bar\fa\bar\fb=\frac{1}{2}\bar\fa+\frac{1}{2}\bar\fb$. 
	This algebra is also baric, being a factor of the baric algebra $\fJ(1)$.
	
	Hence we have the following corollary.
	
	\begin{cor}
		Every $2$-generated algebra $A=\lla a,b\rra$ of Jordan type half, with 
		distinct axes $a$ and $b$, is isomorphic to one of the following:
		\begin{enumerate}[\rm(a)]
			\item $\fJ(\al)$, where $\al=(a,b)$;
			\item the $2B$ algebra $\overline{\fJ(0)}$ (only if $(a,b)=0$);
			\item the baric algebra $\overline{\fJ(1)}$ (only if $(a,b)=1$).
		\end{enumerate} 
		The isomorphism takes $\fa$ and $\fb$ (or $\bar\fa$ and $\bar\fb$) to $a$ and $b$,
		respectively.
	\end{cor}
	
	We next note that, if $\al\neq 1$ then $1_{\fJ}=\frac{2}{\al-1}\bsg$ is the identity
	of the algebra $\fJ$. (Correspondingly, in $\overline{\fJ(0)}$, the 
	identity is $\frac{2}{0-1}\bar\bsg=\bar\fa+\bar\fb$.) 
	
	\medskip
	We will now discuss a different realisation of the majority of these algebras.
	Consider a quadratic space $(V,q)$ over the field $\F$. (Recall that $\F$ is 
	of characteristic not $2$ throughout this paper.) This means that 
	$b(u,v):=\frac{1}{2}(q(u+v)-q(u)-q(v))$ is a symmetric bilinear form on $V$ and 
	we have that $q(u)=b(u,u)$ for all $u\in V$. Define an algebra product on the 
	space $\fS:=\F\oplus V$ as follows: 
	$$(\al+u)(\bt+v):=(\al\bt+b(u,v))+(\al v+\bt u),$$
	for all $\al,\bt\in\F$ and $u,v\in V$.
	
	The resulting algebra, denoted $\fS(V,q)$, is a Jordan algebra and it is 
	known as the \emph{spin factor} algebra or, sometimes, the algebra of 
	\emph{Clifford type}, because it is a subalgebra of the Clifford algebra 
	$Cl(V,q)$ taken with the Jordan product.
	
	Let us now specify the Frobenius form on $\fS=\fS(V,q)$. We take $(u,v)=2b(u,v)$ 
	on $V$ and extend this to the entire $\fS$ by making $1$ orthogonal to $V$ and 
	setting $(1,1)=2$. Then $((\al+u)(\bt+v),\gm+w)=(\al\bt+b(u,v)+\al v+\bt u,
	\gm+w)=2(\al\bt+b(u,v))\gm+2b(\al v+\bt u,w)=2(\al\bt\gm+\gm b(u,v)+
	\bt b(u,w)+\al b(v,w))$. This is clearly symmetric in the three vectors, 
	$\al+u$, $\bt+v$, and $\gm+w$, which means that this is also equal to 
	$(\al+u,(\bt+v)(\gm+w))$, proving that this is a Frobenius form.
	
	All non-zero, non-identity idempotents $a\in\fS$ are of the form 
	$a=\frac{1}{2}(1+u)$, where $u\in V$ satisfies $b(u,u)=1$.\footnote{Note 
		that $v=-u$ also satisfies $b(v,v)=1$, and so 
		$a'=\frac{1}{2}(1+v)=\frac{1}{2}(1-u)$ is also an idempotent.} Furthermore, 
	all these idempotents are primitive and they satisfy $(a,a)=\frac{1}{4}((1,1)+(u,u))=
	\frac{1}{4}(2+2b(u,u))=\frac{1}{4}(2+2)=1$, which is what we usually expect 
	of a Frobenius form.
	
	Finally, let us also mention that $\fS$ is an algebra of Jordan type half if 
	and only if it is generated by primitive idempotents; equivalently, $V$ must 
	be spanned by vectors $u$ with $q(u)=b(u,u)=1$, or equivalently, $(u,u)=2$.
	
	We now again focus on the $2$-generated case, that is, where the algebra is 
	generated by two primitive axes. First take $V=\la u\ra$, where $b(u,u)=1$. 
	Then, setting $a=\frac{1}{2}(1+u)$ and $b=\frac{1}{2}(1-u)$, we have that 
	$\fS(V,q)=\lla a,b\rra\cong\F\oplus\F$, since $ab=0$. This $2B$ algebra 
	is isomorphic to $\overline{\fJ(0)}$. In particular, this is the only algebra 
	we have in the case where $V=\la u\ra$ is $1$-dimensional.
	
	Now take $V=\la u,v\ra$, where $b(u,u)=1=b(v,v)$ and let $\dl=b(u,v)$. In 
	this case, we write $\fS(\dl)$ for $\fS(V,q)$. Let $a=\frac{1}{2}(1+u)$ and 
	$b=\frac{1}{2}(1+v)$. As usual, we use $\al$ to denote $(a,b)$ and note 
	that $\al=(a,b)=(\frac{1}{2}(1+u),\frac{1}{2}(1+v))=\frac{1}{4}(2+2\dl)=
	\frac{1}{2}(1+\dl)$. In particular, $\al=0$ when $\dl=-1$ and $\al=1$ when 
	$\dl=1$.
	
	The following is easy to check.
	
	\begin{prop} \label{spin iso}
		The following hold:
		\begin{enumerate}
			\item[$(a)$] $a$ and $b$ generate $\fS(\dl)$ if and only if $\dl\neq 1$; 
			in particular, in this case, $\fS(\dl)\cong\fJ(\frac{1}{2}(1+\dl))$;
			\item[$(b)$] if $\dl=1$ then $\lla a,b\rra=\la a,b\ra$ is a $2$-dimensional 
			subalgebra isomorphic to $\overline{\fJ(1)}$.
		\end{enumerate} 
	\end{prop} 
	
	The algebra $\fJ(1)$ does not arise in the context of this proposition because 
	$\fS(\dl)$ always has an identity element, while $\fJ(1)$ does not have one.
	
	We note that the isomorphisms in this proposition are again pointed, i.e., 
	they match $a$ and $b$ with $\fa$ and $\fb$ (or $\bar\fa$ and $\bar\fb$), 
	respectively. If we drop pointedness then many of the above algebras become 
	isomorphic even for different values of $\dl$ (equivalently, $\al$). For 
	example, if $\F$ is algebraically closed then all $2$-dimensional non-degenerate 
	orthogonal spaces are isometric, and so all algebras $\fS(\dl)$, $\dl\neq\pm 1$, 
	are isomorphic in this more general sense. (That is, all simple $3$-dimensional 
	algebras of Jordan type half are isomorphic if $\F$ is algebraically closed.)
	
	\medskip
	Next, we construct, under minor additional conditions, a more convenient basis 
	in $\fS=\fS(\dl)$.
	
	\begin{lem} \label{nil basis}
		Suppose that $\F$ contains the roots, $\zeta$ and $\zeta^{-1}$, of the polynomial 
		$x^2-2\dl x+1$ and $\dl\neq\pm 1$. Then 
		\begin{enumerate}
			\item[$(a)$] $\fS(\dl)$ contains exactly two $1$-dimensional nilpotent 
			subalgebras $\la s\ra$ and $\la t\ra$; both of them lie in $V$;
			\item[$(b)$] the vectors $s$, $1$ and $t$ form a basis of $\fS(\dl)$; 
			furthermore, if we scale $s$ and $t$ so that $(s,t)=\frac{1}{4}$ (equivalently,  
			$b(s,t)=\frac{1}{8}$) then the idempotents of $\fS(\dl)$ are $0$, $1$ and all 
			vectors of the form $\xi s+\frac{1}{2}+\xi^{-1}t$ for $\xi\in\F^\sharp=
			\F\setminus\{0\}$.  
		\end{enumerate}
	\end{lem}
	
	\pf First of all, every element of $\fS$ outside of $V$ can be written,
	after scaling, as $1+w$ for $w\in V$. Then $(1+w)^2=1+2w+b(w,w)
	=(1+b(w,w))+2w$. Clearly, this is never zero, and so every $1$-dimensional 
	nilpotent subalgebra must lie in $V$.
	
	Recall from the definition of $\fS$ that $V=\la u,v\ra$, where 
	$b(u,u)=1=b(v,v)$ and $b(u,v)=\dl$. Manifestly, $v^2=b(v,v)=1$ and so 
	$\la v\ra$ is not nilpotent. Every vector of $V$ outside of $\la v\ra$ can be 
	written, after scaling, as $u-\eta v$ for some $\eta\in\F$. Then, $u-\eta v$ 
	generates a $1$-dimensional nilpotent subalgebra if and only if 
	$0=(u-\eta v)^2=b(u-\eta v,u-\eta v)=b(u,u)-2\eta b(u,v)+\eta^2 b(v,v)=
	1-2\dl\eta+\eta^2$. So $\eta$ must be a root of the polynomial 
	$x^2-2\dl x+1$. These roots, $\dl\pm\sqrt{\dl^2-1}$, are distinct if and only 
	if $\dl\neq\pm 1$. This proves (a).
	
	Since $s$ and $t$, being independent, span $V$, we have that $s$, $1$, and 
	$t$ form a basis of $\fS$. Clearly, $(s,t)$ cannot be zero, since 
	$(s,s)=(1,s^2)=(1,0)=0$ and, similarly, $(t,t)=0$. Let us scale 
	$t$ so that $(s,t)=\frac{1}{4}$. Then the multiplication with respect to 
	the basis $s$, $1$, and $t$ is as follows: $s^2=0=t^2$, $st=b(s,t)=\frac{1}{2}(s,t)
	=\frac{1}{8}$ (as $s,t\in V$), and clearly, $1s=s$, $1t=t$, and $1^2=1$.
	
	Turning to idempotents in $\fS$, we first note that every idempotent, apart from 
	$0$ and $1$, has the identity component $\frac{1}{2}$. Taking now an arbitrary 
	element $\xi s+\frac{1}{2}+\lm t$, we have that $(\xi s+\frac{1}{2}+\lm t)^2=
	\xi^2s^2+(\frac{1}{2})^2+\lm^2t^2+2\xi s\frac{1}{2}+2\xi s\lm t
	+2\frac{1}{2}\lm t=0+\frac{1}{4}+0+\xi s+2\xi\lm b(s,t)+\lm t=\frac{1}{4}+
	\xi s+\frac{1}{4}\xi\lm+\lm t=\xi s+\frac{1}{4}(1+\xi\lm)+\lm t$. 
	So $\xi s+\frac{1}{2}+\lm t$ is an idempotent if and only if $\xi\lm=1$, 
	that is, $\lm=\xi^{-1}$, proving (b).\qed
	
	\medskip
	Let us now find in $\fJ=\fJ(\al)$, $\al\neq 0,1$, an equivalent of the basis 
	$\{s,1,t\}$ from $\fS(\dl)$. First of all, $1$ corresponds to 
	$1_{\fJ}=\frac{2}{\al-1}\bsg$. Hence we just need to find in $\fJ$ equivalents 
	of $s$ and $t$. In view of the applications later in the paper, it will be more 
	convenient to do all calculations with the basis $\{\fa,\fb,\fa\fb\}$ of 
	$\fJ$. Since $\fa\fb=\frac{1}{2}\fa+\frac{1}{2}\fb+\bsg$, it is easy to see that 
	we have the following formulae:
	$$\fa(\fa\fb)=\frac{1}{2}(\al\fa+\fa\fb),$$
	$$\fb(\fa\fb)=\frac{1}{2}(\al\fb+\fa\fb),$$
	$$(\fa\fb)^2=\frac{\al}{4}(\fa+\fb+2\fa\fb),$$
	$$(\fa,\fa\fb)=(\fb,\fa\fb)=\al,$$
	and
	$$(\fa\fb,\fa\fb)=\frac{1}{2}\al(\al+1).$$
	
	Let $\zeta,\zeta^{-1}\in\F$ be the roots of the polynomial $x^2-2\dl x+1=x^2-(4\al-2)x+1$. Let us assume for the remainder of this section that $\zeta$ and $\zeta^{-1}$ are elements of $\F$, that is, $x^2-(4\al-2)x+1$ factorises into linear factors over $\F$.
	We note that $\zeta=\zeta^{-1}$ if and only if $\zeta=\pm 1$, and 
	so $4\al-2=\zeta+\zeta^{-1}=\pm 2$, which is only possible when 
	$\al\in\{0,1\}$. So, since we assume that $\al\neq 0,1$, it follows that $\zeta\neq\zeta^{-1}$. 
	
	Define $\mu=-\frac{1+\zeta}{4}$ and $\nu=-\frac{1+\zeta^{-1}}{4}$. Note that 
	$\mu\neq\nu$ since $\zeta\neq\zeta^{-1}$. 
	
	\begin{lem} \label{munu}
		We have the following:
		\begin{enumerate}
			\item[$(a)$] $\mu+\nu=-\al$ and $\mu\nu=\frac{\al}{4}$;
			\item[$(b)$] furthermore, $\mu^2=\frac{\al\zeta}{4}$ and 
			$\nu^2=\frac{\al\zeta^{-1}}{4}$.
		\end{enumerate}
	\end{lem}
	\pf First of all, $\mu+\nu=-\frac{1}{4}(1+\zeta+1+\zeta^{-1})=-\frac{1}{4}(2+4\al-2)
	=-\al$. Also, 
	$\mu\nu=\frac{1}{16}(1+\zeta)(1+\zeta^{-1})=\frac{1}{16}(1+\zeta+\zeta^{-1}+1)=
	\frac{1}{16}(2+4\al-2)=\frac{1}{4}\al$. 
	
	For (b), $\mu^2=\frac{1}{16}(1+2\zeta+\zeta^2)=\frac{1}{16}(1+2\zeta+
	(4\al-2)\zeta-1)=\frac{1}{4}\al\zeta$. Similarly, $\nu^2=\frac{1}{16}(1+
	2\zeta^{-1}+\zeta^{-2})=\frac{1}{16}(1+2\zeta^{-1}+(4\al-2)\zeta^{-1}-1)=
	\frac{1}{4}\al\zeta^{-1}$.\qed 
	
	\medskip
	We now find the two nilpotent subalgebras in $\fJ(\al)$.
	
	\begin{lem} 
		Let $\fd=\mu\fa+\fa\fb+\nu\fb$ and $\fe=\nu\fa+\fa\fb+\mu\fb$. Then
		\begin{enumerate}
			\item[$(a)$] $\fd^2=0=\fe^2$; and
			\item[$(b)$] $(\fd,\fe)=\al(\al-1)^2$.
		\end{enumerate}
	\end{lem}
	
	\pf By symmetry between $\fa$ and $\fb$, it suffices to check that $\fd^2=0$. 
	Using the algebra products above, $\fd^2=(\mu\fa+\fa\fb+\nu\fb)^2=\mu^2\fa+(\fa\fb)^2+\nu^2\fb+2\mu\fa(\fa\fb)+2\nu\fb(\fa\fb)+2\mu\nu\fa\fb=\mu^2\fa+
	\frac{1}{4}\al(\fa+\fb+2\fa\fb)+\nu^2\fb+2\mu\frac{1}{2}(\al\fa+\fa\fb)+
	2\nu\frac{1}{2}(\al\fb+\fa\fb)+2\mu\nu\fa\fb=(\mu^2+\frac{1}{4}\al+\mu\al)\fa+
	(\frac{1}{4}\al+\nu^2+\nu\al)\fb+(\frac{1}{2}\al+\mu+\nu+2\mu\nu)\fa\fb$. 
	Substituting the values from Lemma \ref{munu}, we get $\fd^2=(\frac{1}{4}\al
	\zeta+\frac{1}{4}\al-\frac{1}{4}(1+\zeta)\al)\fa+(\frac{1}{4}\al+\frac{1}{4}
	\al\zeta^{-1}-\frac{1}{4}(1+\zeta^{-1})\al)\fb+(\frac{1}{2}\al-\al+
	2\frac{1}{4}\al)\fa\fb=0x+0y+0\fa\fb=0$.
	
	For part (b), $(\fd,\fe)=(\mu\fa+\fa\fb+\nu\fb,\nu\fa+\fa\fb+\mu\fb)=
	\mu\nu(\fa,\fa)+\mu(\fa,\fa\fb)+\mu^2(\fa,\fb)+\nu(\fa\fb,\fa)+(\fa\fb,\fa\fb)+\mu(\fa\fb,\fb)+\nu^2(\fb,\fa)+\nu(\fb,\fa\fb)+\nu\mu(\fb,\fb)=
	\frac{1}{4}\al+\mu\al+\mu^2\al+\nu\al+\frac{1}{2}(\al^2+\al)+\mu\al+\nu^2\al+
	\nu\al+\frac{1}{4}\al=\al(\frac{1}{4}+\mu+\mu^2+\nu+\frac{1}{2}\al+\frac{1}{2}+
	\mu+\nu^2+\nu+\frac{1}{4})=\al(1+2(\mu+\nu)+\mu^2+\nu^2+\frac{1}{2}\al)=
	\al(1-2\al+\frac{1}{4}\al\zeta+\frac{1}{4}\al\zeta^{-1}+\frac{1}{2}\al)=
	\al(1-2\al+\frac{1}{4}\al(4\al-2)+\frac{1}{2}\al)=\al(\al-1)^2$.\qed
	
	\medskip
	Clearly, since $\mu\neq\nu$, the elements $\fd$ and $\fe$ are linearly independent. 
	So the two $1$-dimensional nilpotent subalgebras in $\fJ$ are $\la\fd\ra$ and 
	$\la\fe\ra$. In particular, we can choose elements $\fs$ and $\ft$, as in Lemma 
	\ref{nil basis}, by scaling $\fd$ and $\fe$. Namely, we set $\fs=\frac{1}{4(\al-1)}\fd$ 
	and $\ft=\frac{1}{\al(\al-1)}\fe$. Then we have that $(\fs,\ft)=\frac{1}{4}$ and hence 
	also $\fs\ft=\frac{1}{8}1_{\fJ}=\frac{2}{8(\al-1)}\bsg=\frac{1}{4(\al-1)}\bsg$.
	
	Let us now see what $\fa$ and $\fb$ look like in terms of the basis $\fs$, $1_{\fJ}$, 
	and $\ft$.
	
	\begin{lem}
		We have that
		\begin{enumerate}
			\item[$(a)$] $\fa=\nu^{-1}\fs+\frac{1}{2}1_{\fJ(\al)}+\nu\ft$; and
			\item[$(b)$] $\fb=\mu^{-1}\fs+\frac{1}{2}1_{\fJ(\al)}+\mu\ft$.
		\end{enumerate}
	\end{lem}
	
	\pf We let $\fa=\xi\fs+\frac{1}{2}1_{\fJ}+\xi^{-1}\ft$ for some $0\neq\xi\in\F$. Then, 
	on the one hand, $(\fa,\ft)=\frac{1}{4}\xi+0+0=\frac{1}{4}\xi$, that is, $\xi=4(\fa,\ft)$.
	On the other hand,
	\begin{align*}
		4(\fa,\ft)&=4\left(\fa,\frac{1}{\al(\al-1)}(\nu\fa+\fa\fb+\mu\fb)\right)=4\frac{\nu+\al+\mu\al}
		{\al(\al-1)}
		\\&=4\frac{\nu-\nu-\mu+\mu\al}{\al(\al-1)}=4\frac{\mu(\al-1)}{\al(\al-1)}=4\frac{\mu}{\al}=
		\nu^{-1}.
	\end{align*}
	
	Thus, $\xi=\nu^{-1}$ and $\fs=\nu^{-1}\fs+\frac{1}{2}1_{\fJ}+\nu\ft$. The calculation for $\fb$ 
	is completely similar.\qed
	
	\medskip
	We now look at the action of the Miyamoto involutions of axes from $\fJ$.
	
	\begin{lem} \label{action on J}
		If $\fd=\xi\fs+\frac{1}{2}1_{\fJ}+\xi^{-1}\ft$ is an axis in $\fJ(\al)$ for some $0\neq\xi\in\F$ then 
		\begin{enumerate}
			\item[$(a)$] $\fs^{\tau_{\fd}}=\xi^{-2}\ft$; and 
			\item[$(b)$] $\ft^{\tau_{\fd}}=\xi^2\fs$.
		\end{enumerate}
	\end{lem}
	
	\pf Note that $(\fd,\fs)=(\xi\fs+\frac{1}{2}1_{\fJ}+\xi^{-1}\ft,\fs)=\frac{1}{4}\xi^{-1}$ and similarly 
	$(\fd,\ft)=\frac{1}{4}\xi$. Furthermore, $\fs$ and $\ft$ are the only elements of $\la\fs\ra$ and $\la\ft\ra$, 
	respectively, having these values of Frobenius form with $\fd$. If $\tau_{\fd}$ preserves 
	the nilpotent $1$-spaces $\la\fs\ra$ and $\la\ft\ra$ then $\tau_{\fd}$ must fix both $\fs$ and $\ft$, 
	as it preserves the Frobenius form. Since $\tau_{\fd}$ also fixes the identity $1_{\fJ}$, we conclude that 
	$\tau_{\fd}$ acts as identity on $\fJ$, which is clearly not true. (Since $\fJ$ is $3$-dimensional, 
	$\fJ_{\frac{1}{2}}(\fd)\neq 0$.)
	
	We conclude that $\tau_{\fd}$ switches $\la\fs\ra$ and $\la\ft\ra$. Hence, there exists an $\om\in\F$ such 
	that $\fs^{\tau_{\fd}}=\om\ft$. Furthermore, $(\fd,\om\ft)=(\fd,\fs)=\frac{1}{4}\xi^{-1}$. Since 
	$(\fd,\ft)=\frac{1}{4}\xi$, we conclude that 
	$\om=\xi^{-2}$. That is, $\fs^{\tau_{\fd}}=\xi^{-2}\ft$. Similarly, we get that $\ft^{\tau_{\fd}}
	=\xi^2\fs$.\qed
	
	\medskip
	Now we set $\rho=\tau_{\fa}\tau_{\fb}$. As a consequence of the above, we deduce the action of $\rho$ 
	on $\fJ$. 
	
	\begin{cor} \label{on R}
		The following hold:
		\begin{enumerate}
			\item[$(a)$] $\rho$ acts on $\la\fs\ra$ as the scalar $\zeta^{-2}$ and, symmetrically, $\rho$ 
			acts on $\la\ft\ra$ as the scalar $\zeta^2$;
			\item[$(b)$] $\rho$ acts on $\la 1_{\fJ(\al)}\ra$ as the identity.
		\end{enumerate}
	\end{cor}
	
	\pf Since $\fa=\nu^{-1}\fs+\frac{1}{2}1_{\fJ}+\nu\ft$ and $\fb=\mu^{-1}\fs+\frac{1}{2}1_{\fJ}+\mu\ft$, 
	we have that $\fs^\rho=(\fs^{\tau_{\fa}})^{\tau_{\fb}}=(\nu^2\ft)^{\tau_{\fb}}=\mu^{-2}\nu^2 \fs$. By 
	Lemma \ref{munu} (b), $\mu^{-2}\nu^2 =\zeta^{-2}$, as claimed. Similarly, $\ft^\rho=
	(\ft^{\tau_{\fa}})^{\tau_{\fb}}=\nu^{-2}\mu^2\ft=\zeta^2\ft$.
	
	Part (b) is clear as both $\tau_{\fa}$ and $\tau_{\fb}$ fix $1_{\fJ}$.\qed

	\section{Three axes} \label{3-gen}
	
	Next we need to review the properties of algebras of Jordan type half that 
	are generated by three axes. Here we follow \cite{gs}. 
	
	Let $\fT=\fT(\al,\bt,\gm,\psi)$, $\al,\bt,\gm,\psi\in\F$, denote the commutative 
	algebra over $\F$ generated by idempotents $\fa$, $\fb$, and $\fc$, and having basis 
	$\fa$, $\fb$, $\fc$, $\fa\fb$, $\fa\fc$, $\fb\fc$, $\fa(\fb\fc)$, $\fb(\fa\fc)$, $\fc(\fa\fb)$. To specify the 
	algebra product on $\fT$, it will be convenient to use the following notation: 
	$\al=(\fa,\fb)=(\fb,\fa)$, $\bt=(\fb,\fc)=(\fc,\fb)$, $\gm=(\fc,\fa)=(\fa,\fc)$, and $\psi=(\fa,\fb\fc)
	=(\fb,\fc\fa)=(\fc,\fa\fb)$. (These are indeed the values of the Frobenius form.)
	
	Using this notation, multiplication on the above basis, where we skip the 
	obvious products, is given in Table \ref{t:mult}, where $\{u,v,w\}=\{\fa,\fb,\fc\}$.
	\begin{table}[t!] \label{t:mult}
		\renewcommand{\arraystretch}{1.7}
		\caption{Multiplication in $\fT(\al,\bt,\gm,\psi)$}
		\medskip
		\begin{tabular}{|l|}
			\hline
			$u(uv)=\frac{1}{2}((u,v)u+uv)$\\
			\hline
			$u(u(vw))=\frac{1}{2}((u,vw)u+u(vw))$\\
			\hline
			$u(v(uw))=\frac{1}{4}((u,vw)u+(u,w)uv+(u,v)uw+u(vw)+v(uw)-w(uv))$\\
			\hline
			$(uv)^2=\frac{1}{4}(u,v)(u+v+2uv)$\\
			\hline
			$(uv)(vw)=\frac{1}{4}((v,uw)v+(v,w)uv+(u,v)vw+u(vw)-v(uw)+w(uv))$\\
			\hline
			$(uv)(u(vw))=\frac{1}{8}(((u,v)(u,w)+(u,vw))u+(u,vw)v+2(u,vw)uv$\\
			\quad $+(u,v)vw+2(u,v)u(vw))$\\
			\hline
			$(uv)(w(uv))=\frac{1}{8}( (u,v)(v,w){u}+(u,v)(u,w){v}+4(uv,w){uv}$\\ 
			\quad $+(u,v){vw}+(u,v){uw}-2(u,v){u(vw)}-2(u,v){v(uw)}+4(u,v){w(uv)} )$\\
			\hline
			$(u(vw))^2=\frac{1}{16}((v,w)((u,v)+(u,w)+2(u,vw))u+(v,w)(u,w)v$\\
			\quad $+(u,v)(v,w)w+4(u,vw)vw+(2(v,w)+8(u,vw))u(vw)$\\
			\quad $-2(v,w)v(uw)-2(v,w)w(uv))$\\
			\hline
			$(u(vw))(v(uw))=\frac{1}{16}((v,w)((u,w)+(u,vw))u+(u,w)((v,w)+(u,vw))v$\\
			\quad $+(u,v)(u,vw){w}+2(u,vw){uv}+2((u,v)(u,w)+(u,vw))vw$\\
			\quad $+2((u,v)(v,w)+(u,vw))uw+(4(u,vw)-(u,v)+(v,w)$\\
			\quad $-(u,w))u(vw)+(4(u,vw)-(u,v)-(v,w)+(u,w))v(uw)$\\
			\quad $+((u,v)-(v,w)-(u,w)-4(u,vw))w(uv))$\\
			\hline 
		\end{tabular}
	\end{table}
	This algebra $\fT$ is of Jordan type half, in fact, a Jordan algebra, for all 
	choices of the free parameters $\al$, $\bt$, $\gm$, and $\psi$. Furthermore, 
	$\fa$, $\fb$, and $\fc$ are primitive axes. The unique Frobenius form on $\fT$ has 
	values shown in Table \ref{t:gram}, where again $\{u,v,w\}=\{\fa,\fb,\fc\}$. All 
	expressions are in terms of the basic values $\al$, $\bt$, $\gm$, and $\psi$.    
	\begin{table}[t!] \label{t:gram}
		\renewcommand{\arraystretch}{1.7}
		\caption{Frobenius form on $\fT(\al,\bt,\gm,\psi)$}
		
		\medskip
		\begin{tabular}{|l|}
			\hline
			$(u,uv)=(u,v)$\\
			\hline
			$(u,u(vw))=(u,vw)$\\
			\hline
			$(u,v(uw))=\frac{1}{2}((u,v)(u,w)+(u,vw))$\\
			\hline
			$(uv,uv)=\frac{1}{2}(u,v)((u,v)+1)$\\
			\hline
			$(uv,vw)=\frac{1}{2}((u,v)(v,w)+(u,vw))$\\
			\hline
			$(uv,u(vw))=\frac{1}{4}((u,v)(v,w)+(u,vw)+2(u,v)(u,vw))$\\
			\hline
			$(uv,w(uv))=\frac{1}{4}(u,v)((v,w)+(u,w)+2(u,vw))$\\
			\hline
			$(u(vw),u(vw))=\frac{1}{8}((u,v)(v,w)+(u,w)(v,w)+2(b,c)(u,bc)+4(u,vw)^2)$\\
			\hline
			$(u(vw),v(uw))=\frac{1}{8}(2(u,v)(u,w)(v,w)+(u,v)(u,vw)+(u,w)(u,vw)$\\
			\quad $+(v,w)(u,vw)+(u,w)(v,w)+2(u,vw)^2)$\\
			\hline
		\end{tabular}
	\end{table}
	
	\medskip
	The following is the main result, Proposition 5, from \cite{gs}.
	
	\begin{thm} \label{3-gen main}
		Suppose that $A$ is an algebra of Jordan type half over $\F$ generated by three 
		primitive axes $a$, $b$, and $c$. Let $\al=(a,b)$, $\bt=(b,c)$, $\gm=(c,a)$, and 
		$\psi=(a,bc)$. Then there exists a (unique) surjective homomorphism from
		$\fT(\al,\bt,\gm,\psi)$ onto $A$ sending $\fa$, $\fb$, and $\fc$ to $a$, $b$, and $c$, 
		respectively.
	\end{thm}
	
	In other words, $\fT$ is the universal $3$-generated algebra 
	of Jordan type half.\footnote{In fact, we obtain the universal algebra when we 
		define $\fT$ over the polynomial ring $\Z[\frac{1}{2}][\al,\bt,\gm,\psi]$, where 
		$\al$, $\bt$, $\gm$, and $\psi$ are indeterminates. Then the homomorphism 
		from $\fT$ to $A$ is accompanied by the evaluation map $\Z[\frac{1}{2}][\al,\bt,
		\gm,\psi]\to\F$ sending the indeterminates to the suitable values of the Frobenius 
		form. Recall that our field $\F$ cannot be of characteristic $2$, and so the image 
		of $\frac{1}{2}$ is defined.}
	
	We note that the subalgebra $\fJ=\lla\fa,\fb\rra$ of $\fT$ is a copy of $\fJ(\al)$ 
	although here the multiplication is written with respect to the basis 
	$\fa$, $\fb$, and $\fa\fb$, instead of $\fa$, $\fb$, and $\bsg=\fa\fb-
	\frac{1}{2}(\fa+\fb)$. For the remainder of this section, let us assume that 
	$\al\neq 1,0$. Then, by Theorem \ref{all ideals}, $\fJ$ is simple of dimension $3$. 
	Furthermore, as in the preceding section, $\fJ$ contains elements 
	$\fs$ and $\ft$, such that $\fs^2=0=\ft^2$ and $\fa=\nu^{-1}\fs+\frac{1}{2}1_{\fJ}+
	\nu\ft$ and 
	$\fb=\mu^{-1}\fs+\frac{1}{2}1_{\fJ}+\mu\ft$, where $\mu=-\frac{1+\zeta}{4}$, 
	$\nu=-\frac{1+\zeta^{-1}}{4}$ and $\zeta$, $\zeta^{-1}$ are the roots of 
	$x^2-(4\al-2)x+1$. (We again assume that $\zeta$ and $\zeta^{-1}$ are elements of 
	$\F$.)
	
	In Lemma \ref{action on J}, we described the action of $\rho=\tau_{\fa}\tau_{\fb}$ 
	on $\fJ=\la\fa,\fb\ra$. Next we consider how this action extends to the entire $\fT$. 
	Let $\fW=\fJ^\perp$ be the orthogonal complement of $\fJ$ in $\fT$. Naturally, the 
	complement is taken with respect to the Frobenius form on $\fT$. Since $\fJ$ is 
	simple, the Frobenius form on it is non-degenerate, which means that $\fW\cap\fJ=0$, 
	and hence
	$$\fT=\fJ\oplus\fW.$$
	
	Let us establish the basic properties of this decomposition. For $u\in\fT$, we will 
	denote by $u^\circ$ the projection of $u$ to $\fW$.
	
	\begin{lem} \label{module}
		The complement $\fW$ is a $\fJ$-module, i.e., $\fW\fJ\subseteq\fW$. Furthermore, 
		the projection map $u\mapsto u^\circ$ commutes with the action of $\fJ$, that is, 
		for $u\in\fT$ and $v\in\fJ$, we have that $(uv)^\circ=u^\circ v$.
	\end{lem}
	
	\pf For $u\in\fW$ and $v,w\in\fJ$, we have that $(uv,w)=(u,vw)=0$, since $vw\in\fJ$. 
	Since this is true for all $w$, we have that $uv\in\fW$, i.e., $\fW$ is indeed 
	a $\fJ$-module.
	
	For the second claim, let $u\in\fT$ and $v\in\fJ$. Then $u=w+u^\circ$ for some 
	$w\in\fJ$. Hence $uv=(w+u^\circ)v=wv+u^\circ v$. Clearly, the first summand is in $\fJ$ 
	and the second summand is in $\fW$. Thus, $(uv)^\circ$ is indeed equal to $u^\circ 
	v$.\qed
	
	\medskip
	First we look at the adjoint action of $\fs$ and $\ft$ on $\fW$. Let 
	$\om=\ad_{\fs}|_{\fW}$ and $\chi=\ad_{\ft}|_{\fW}$ be the restrictions of these 
	adjoint maps to $\fW$.
	
	Note that the projection is a surjective map from $\fT$ onto $\fW$. In particular, 
	the image of a basis of $\fT$, such as $\{\fa,\fb,\fc,\fa\fb,\fb\fc,\fa\fc,\fa(\fb\fc), 
	\fb(\fa\fc),\fc(\fa\fb)\}$, spans $\fW$. Clearly, $\fa^\circ=\fb^\circ=(\fa\fb)^\circ=0$.
	Hence $\fc^\circ$, $(\fb\fc)^\circ$, $(\fa\fc)^\circ$, $(\fa(\fb\fc))^\circ$, 
	$(\fb(\fa\fc))^\circ$, and $(\fc(\fa\fb))^\circ$ span $\fW$. On the other hand, 
	$\dim\fW=\dim\fT-\dim\fJ=6$, and so these six elements form a basis of $\fW$. 
	
	\begin{lem} \label{dimensions}
		The following hold:
		\begin{enumerate}
			\item[$(a)$] $\im\om$ and $\im\chi$ are of dimension $2$;
			\item[$(b)$] $\im\om\subseteq\ker\om$ and, symmetrically, $\im\chi\subseteq\ker\chi$;
			\item[$(c)$] $\fD:=\ker\om\cap\ker\chi$ is of dimension $2$.
		\end{enumerate}
	\end{lem} 
	
	\pf The action of $\om$ with respect to the above basis of $\fW$ is given by 
	the matrix
	$$
	S=\frac{1}{4(\al-1)}\left(
	\renewcommand\arraycolsep{3pt}
	\begin{array}{cccccc}
		0&\nu&\mu&0&0&1\\ 
		0&-\frac{1}{4}(\mu-\nu)&0&-\frac{1}{4}\zeta&-\frac{1}{4}&\frac{1}{4}\\
		0&0&\frac{1}{4}(\mu-\nu)&-\frac{1}{4}&-\frac{1}{4}\zeta^{-1}& \frac{1}{4}\\
		0&\frac{1}{4}\mu\nu(1-\zeta^{-1})&0&\frac{1}{4}\mu&\frac{1}{4}\nu&-\frac{1}{4}\nu\\
		0&0&\frac{1}{4}\mu\nu(1-\zeta)&\frac{1}{4}\mu&\frac{1}{4}\nu&-\frac{1}{4}\mu\\
		0&\frac{1}{4}\mu\nu(1-\zeta^{-1})&\frac{1}{4}\mu\nu(1-\zeta)&\frac{1}{2}\mu&\frac{1}{2}\nu&
		-\frac{1}{4}(\mu+\nu) 
	\end{array}
	\right)
	$$
	and, symmetrically, the action of $\chi$ is given by the matrix 
	$$
	T=\frac{1}{\al(\al-1)}\left(
	\renewcommand\arraycolsep{3pt}
	\begin{array}{cccccc}
		0&\mu&\nu&0&0&1\\ 
		0&\frac{1}{4}(\mu-\nu)&0&-\frac{1}{4}\zeta^{-1}&-\frac{1}{4}&\frac{1}{4}\\
		0&0&\frac{1}{4}(\nu-\mu)&-\frac{1}{4}&-\frac{1}{4}\zeta&\frac{1}{4}\\
		0&\frac{1}{4}\mu\nu(1-\zeta)&0&\frac{1}{4}\nu&\frac{1}{4}\mu&-\frac{1}{4}\mu \\
		0&0&\frac{1}{4}\mu\nu(1-\zeta^{-1})&\frac{1}{4}\nu&\frac{1}{4}\mu&-\frac{1}{4}\nu \\
		0&\frac{1}{4}\mu\nu(1-\zeta)&\frac{1}{4}\mu\nu(1-\zeta^{-1})&\frac{1}{2}\nu&\frac{1}{2}\mu&
		-\frac{1}{4}(\mu+\nu) 
	\end{array}
	\right).
	$$
	
	This can be computed by hand using Lemma \ref{module} and the multiplication rules in 
	Table \ref{t:mult}, since $\fs=\frac{1}{4(\al-1)}(\mu\fa+\fa\fb+\nu\fb)$ and $\ft=
	\frac{1}{\al(\al-1)}(\nu\fa+\fa\fb+\mu\fb)$. However, we found $S$ and $T$ in GAP in 
	order to avoid likely computational errors. Given the matrices above, we can see that 
	the rank of both $S$ and $T$ is $2$ regardless of the values of $\zeta$, $\mu$ and $\nu$. 
	Indeed, looking at the final two entries, we see that, say, the first two rows of $S$ 
	are linearly independent and we can also easily express rows 3--6 as linear combinations 
	of the first two rows. This proves part (a). Also, it is easy to check that $S^2=0=T^2$, 
	which means that (b) also holds.
	
	Furthermore, this means that both $\ker\om$ and $\ker\chi$ are $4$-dimensional, and so 
	$\fD=\ker\om\cap\ker\chi$ is at least $2$-dimensional. On the other hand, the matrix $ST$ 
	can be computed and it contains a $2\times 2$ submatrix (rows 2 and 3, columns 2 and 3)
	$$
	\frac{\al}{2^{14}(\al-1)^2}\left(
	\begin{array}{cc}
		\zeta^{-1}(\zeta-1)^3&0\\
		0&\zeta(\zeta^{-1}-1)^3\\
	\end{array}
	\right),
	$$
	where we simplified the expressions in the matrix using Lemma \ref{munu}. The determinant 
	of this submatrix is $-\frac{\al(\zeta-1)^6}{2^{14}(\al-1)^2\zeta^3}$ and it is not zero. 
	Indeed, $\al\neq 0,1$ by assumption and $\zeta\neq 1$, since we know that 
	$\zeta\neq\zeta^{-1}$.
	
	This shows that $ST$ is of rank at least $2$ and, therefore, $\im\om\cap\ker\chi=0$. 
	Symmetrically, we also have that $\im\chi\cap\ker\om=0$. This shows that 
	$\fD=\ker\om\cap\ker\chi$ has dimension exactly $2$, proving (c).\qed
	
	\medskip
	Our next goal is to find the eigenvalues of $\rho$ acting on $\fW$. (Note that $\rho$ 
	leaves $\fJ$ invariant, and so it also leaves its orthogonal complement, $\fW$, invariant.) 
	Using the previous lemma, we can further decompose $\fW$. Let $\fC=\im\om+\im\chi$. Then we 
	have that $\fC$ is of dimension $4$ and $\fW=\fC\oplus\fD$. 
	
	\begin{lem}
		Both $\fC$ and $\fD$ are $\fJ$-modules, i.e., $\fC\fJ\subseteq\fC$ and $\fD\fJ\subseteq\fD$.
	\end{lem}
	
	\pf Note that $1_{\fJ}$ associates with $\fs$ and $\ft$ by the Seress Lemma\footnote{Note that 
		$\fT$ 
		is a Jordan algebra and hence every non-zero idempotent in it is an axis. Hence the Seress 
		Lemma (Proposition \ref{Seress Lemma}) applies to $1_{\fJ}$.}. Since $1_{\fJ}$ associates 
	with $\fs$ and $\ft$, the linear map $\ad_{1_{\fJ}}|_{\fW}$ commutes with both $\om$ and 
	$\chi$, and hence $\im\om$, $\im\chi$, and $\fD=\ker\om\cap\ker\chi$ are all invariant under 
	$\ad_{1_{\fJ}}|_{\fW}$.
	
	On the other hand, both $\fC$ and $\fD$ are clearly invariant under $\om$ and $\chi$, 
	and so $\fC$ and $\fD$ are $\fJ$-modules, since $\{\fs,1_{\fJ},\ft\}$ is a basis of 
	$\fJ$.\qed
	
	\medskip
	We can now determine the action of $\rho$ on $\fD$.
	
	\begin{lem} \label{on D}
		The map $\rho=\tau_{\fa}\tau_{\fb}$ acts on $\fD$ as identity.
	\end{lem}
	
	\pf Recall that both $\fs$ and $\ft$ act on $\fD$ trivially, which implies that 
	$\fa=\nu^{-1}\fs+\frac{1}{2}1_{\fJ}+\nu\ft$ and $\fb=\mu^{-1}\fs+\frac{1}{2}1_{\fJ}+\mu\ft$ 
	act on $\fD$ in the same way, which implies that $\tau_{\fa}|_{\fD}=\tau_{\fb}|_{\fD}$. Thus, 
	$\rho|_{\fD}=\tau_{\fa}|_{\fD}\tau_{\fb}|_{\fD}=
	\tau_{\fa}^2|_{\fD}$ is the identity map.\qed
	
	\medskip
	Next we focus on $\fC$. Let $w\in\im\om$. For generality, let $u=\xi\fs+\frac{1}{2}1_{\fJ}+
	\xi^{-1}\ft$, where $\xi\in\F$, be a primitive axis of $\fT$ contained in $\fJ$. (Say, $u=\fa$ 
	or $\fb$.) We note that $\fW$ decomposes as $\fW_0(u)\oplus\fW_{\frac{1}{2}}(u)$ by primitivity, 
	since $\fW_1(u)=\la u\ra\subseteq\fJ$. In particular, $w=w_0+w_{\frac{1}{2}}$ with 
	$w_0\in\fW_0(u)$ and $w_{\frac{1}{2}}\in\fW_{\frac{1}{2}}(u)$.
	
	\begin{lem} \label{s and t}
		We have the following:
		\begin{enumerate}
			\item[$(a)$] $\ad_{1_{\fJ}}$ acts on $\fC$ as the scalar $\frac{1}{2}$; and
			\item[$(b)$] $\chi\om$ acts on $\im\om$ as the scalar $\frac{1}{16}$ and, symmetrically, 
			$\om\chi$ acts as $\frac{1}{16}$ on $\im\chi$.
		\end{enumerate}
	\end{lem}
	
	\pf We compute: $\frac{1}{2}w_{\frac{1}{2}}=wu=w(\xi\fs+\frac{1}{2}1_{\fJ}+\xi^{-1}\ft)=\xi w\fs+
	\frac{1}{2}w1_{\fJ}+\xi^{-1}w\ft$. Since $w\in\im\om\subseteq\ker\om$, we have that $\xi w\fs=0$. 
	Thus, $w_{\frac{1}{2}}=w1_{\fJ}+2\xi^{-1}w\ft$. Also, $\frac{1}{2}w1_{\fJ}+\xi^{-1}w\ft=
	\frac{1}{2}w_{\frac{1}{2}}=w_{\frac{1}{2}}u=(w1_{\fJ}+2\xi^{-1}w\ft)(\xi\fs+\frac{1}{2}1_{\fJ}+
	\xi^{-1}\ft)=\xi(w1_{\fJ})\fs+2(w\ft)\fs+\frac{1}{2}(w1_{\fJ})1_{\fJ}+\xi^{-1}(w\ft)1_{\fJ}+
	\xi^{-1}(w1_{\fJ})\ft+2\xi^{-2}(w\ft)\ft$. We now notice that the first and last terms in this 
	sum are zero. Indeed, $w1_{\fJ}\in\im\om\subseteq\ker\om$ and $w\ft\in\im\chi\subseteq\ker\chi$. 
	Furthermore, since $1_{\fJ}$ associates with $\ft$, the terms $\xi^{-1}(w\ft)1_{\fJ}$ and 
	$\xi^{-1}(w1_{\fJ})\ft$ are equal. 
	Thus, 
	$$\frac{1}{2}w1_{\fJ}+\xi^{-1}w\ft=2(w\ft)\fs+\frac{1}{2}(w1_{\fJ})1_{\fJ}+2\xi^{-1}(w\ft)1_{\fJ}.$$
	
	Note that $\frac{1}{2}w1_{\fJ}$, $2(w\ft)\fs$, and $\frac{1}{2}(w1_{\fJ})1_{\fJ}$ are in $\im\om$, 
	while $\xi^{-1}w\ft$ and $2\xi^{-1}(w\ft)1_{\fJ}$ are in $\im\chi$. Since $\im\chi\cap\im\om\subseteq 
	\im\chi\cap\ker\om=0$, we deduce two shorter equalities:

	$$\frac{1}{2}w1_{\fJ}=2(w\ft)\fs+\frac{1}{2}(w1_{\fJ})1_{\fJ}$$
	and
	$$\xi^{-1}w\ft=2\xi^{-1}(w\ft)1_{\fJ}.$$
	
	In particular, the second equality implies (a). Indeed, cancelling $\xi^{-1}$ and dividing 
	by $2$, we get that $\frac{1}{2}(w\ft)=(w\ft)1_{\fJ}$. As $w$ runs through $\im\om$, the product 
	$w\ft$ runs through $\im\chi$, which means that $\ad_{1_{\fJ}}$ acts on $\im\chi$ as the scalar 
	$\frac{1}{2}$. By symmetry, the same is true for $\im\om$, and so indeed (a) holds.
	
	In turn, (a) allows to simplify the first equality above. Namely, $\frac{1}{2}w1_{\fJ}=
	\frac{1}{4}w$ and $\frac{1}{2}(w1_{\fJ})1_{\fJ}=\frac{1}{8}w$. This yields the equality 
	$$\frac{1}{8}w=2(w\ft)\fs,$$
	that is, $\chi\om$ acts as $\frac{1}{16}$ on $\im\om$. Symmetrically, $\om\chi$ acts 
	as $\frac{1}{16}$ on $\im\chi$, and so (b) holds.\qed
	
	\medskip
	Note that $\om\chi$ acts as zero on $\im\om$ and $\chi\om$ acts as zero on $\im\chi$.
	
	Before we start computing the action of $\rho$ on $\fC$, we note that, by Lemma \ref{action on J},
	$\tau_{\fa}$ and $\tau_{\fb}$ switch $\la\fs\ra$ and $\la\ft\ra$, which implies that they 
	also switch $\im\om$ and $\im\chi$.
	
	\begin{lem} \label{action of tau}
		Let again $u=\xi\fs+\frac{1}{2}1_{\fJ}+\xi^{-1}\ft\in\fJ$ be a primitive axis of $\fT$.
		\begin{enumerate}
			\item[$(a)$] If $w\in\im\om$ then $w^{\tau_u}=-4\xi^{-1}w\ft$; and symmetrically, 
			\item[$(b)$] if $w\in\im\chi$ then $w^{\tau_u}=-4\xi w\fs$.
		\end{enumerate}
	\end{lem}
	
	\pf Let $w\in\im\om$. We have earlier established the equality 
	$$w_{\frac{1}{2}}=w1_{\fJ}+2\xi^{-1}w\ft=\frac{1}{2}w+2\xi^{-1}w\ft.$$
	Clearly, $\tau_u$ negates $w_{\frac{1}{2}}$. Hence, we deduce 
	$$\frac{1}{2}w^{\tau_u}+2\xi^{-1}(w\ft)^{\tau_u}=w_{\frac{1}{2}}^{\tau_u}=-w_{\frac{1}{2}}=
	-\frac{1}{2}w-2\xi^{-1}w\ft.$$
	
	Since $\tau_u$ switches $\im\om$ and $\im\chi$, we have 
	$$\frac{1}{2}w^{\tau_u}=-2\xi^{-1}w\ft$$
	and
	$$2\xi^{-1}(w\ft)^{\tau_u}=-\frac{1}{2}w.$$
	
	The former clearly yields (a), while the second yields (b), since $w\ft$ is an arbitrary 
	element of $\im\chi$ and $(w\ft)^{\tau_u}=-\frac{1}{4}\xi w=-\frac{1}{4}\xi[16(w\ft)\fs]=
	-4\xi(w\ft)\fs$ by Lemma \ref{s and t}.\qed
	
	\medskip
	Finally, we are ready to establish the action of $\rho=\tau_{\fa}\tau_{\fb}$ on $\im\om$ and 
	$\im\chi$.
	
	\begin{lem} \label{on C}
		The map $\rho$ acts on $\im\om$ as the scalar $\zeta^{-1}$ and on $\im\chi$ as $\zeta$.
	\end{lem}
	
	\pf Take $w\in\im\om$. Applying Lemma \ref{action of tau} with $u=\fa$, we get that 
	$w^{\tau_{\fa}}=-4\nu w\ft$, since $\xi=\nu^{-1}$ for the axis $\fa$. Now applying Lemma 
	\ref{action of tau} with $u=\fb$, we get $w^\rho=(w^{\tau_{\fa}})^{\tau_{\fb}}=
	(-4\nu w\ft)^{\tau_{\fb}}=
	-4\nu(-4\mu^{-1}(w\ft)\fs)=16\nu\mu^{-1}\frac{1}{16}w=\nu\mu^{-1}w$. It remains to notice 
	that $\nu\mu^{-1}=\frac{\nu^2}{\mu\nu}=\frac{\frac{\al\zeta^{-1}}{4}}{\frac{\al}{4}}=
	\zeta^{-1}$ using Lemma \ref{munu}. This proves the first claim and the second claim now 
	follows by symmetry.\qed
	
	\medskip
	Now that we understand what happens in the universal $9$-dimensional algebra 
	$\fT(\al,\bt,\gm,\psi)$, we can deal with the general $3$-generated case. Suppose 
	that $A=\lla a,b,c\rra$ is a $3$-generated algebra of Jordan type half. We assume 
	that $\al=(a,b)\neq 0,1$. Set $\bt=(b,c)$, $\gm=(a,c)$, $\psi=(a,bc)$ and let 
	$\fT=\fT(\al,\bt,\gm,\psi)$.
	
	By Theorem \ref{3-gen main}, there exists a surjective homomorphism $\phi:\fT\to A$ 
	such that $\phi(\fa)=a$, $\phi(\fb)=b$, and $\phi(\fc)=c$. Clearly, $\phi$ maps 
	$\fJ=\lla\fa,\fb\rra$ isomorphically onto $J=\lla a,b\rra$, since $\al\neq 0,1$ 
	and so both $\fJ$ and $J$ are $3$-dimensional simple. Set $s=\phi(\fs)$ and 
	$t=\phi(\ft)$, where $\fs$ and $\ft$ are the nilpotent elements of $\fJ$ defined earlier.
	
	Let $W=J^\perp$ be the orthogonal complement of $J$ in $A$. Since $J$ is simple, we 
	have again that $J\cap W=0$, and so $A=J\oplus W$, with $W$ being a $J$-module. Also, 
	we have that $\phi$ maps $\fW$ onto $W$.
	
	By a slight abuse of notation, let $\om=\ad_s|_W$, $\chi=\ad_t|_W$, and $\rho=\tau_a\tau_b$. 
	Clearly, $\phi$ transforms the maps $\om$, $\chi$, and $\rho$ from $\fT$ into their 
	namesakes in $A$. By Lemma \ref{on R} applied to $J$, $\rho$ acts as $\zeta^{-2}$ on 
	$\la s\ra$, as $1$ on $\la 1_J\ra$, and as $\zeta^2$ on $\la t\ra$.
	
	\begin{prop} \label{3-gen W}
		In $A=\lla a,b,c\rra$, we have the following:
		\begin{enumerate}
			\item[$(a)$] $\chi\om$ acts as $\frac{1}{16}$ on $\im\om$ and, symmetrically, $\om\chi$ 
			acts as $\frac{1}{16}$ on $\im\chi$;
			\item[$(b)$] $W=\im\om\oplus(\ker\om\cap\ker\chi)\oplus\im\chi$;
			\item[$(c)$] furthermore, $\dim(\im\om)=\dim(\im\chi)\leq 2$ and $\dim(\ker\om\cap\ker\chi)
			\leq 2$;
			\item[$(d)$] finally, $\rho$ acts as $\zeta^{-1}$ on $\im\om$, as $1$ on $\ker\om\cap\ker\chi$, 
			and as $\zeta$ on $\im\chi$.
		\end{enumerate}
	\end{prop}
	
	\pf Part (a) follows from Lemma \ref{s and t}(b). In $\fT$, we have that 
	$\fW=\im\om\oplus(\ker\om\cap\ker\chi)\oplus\im\chi$ and, by Lemmas \ref{on D} 
	and \ref{on C}, $\rho$ acts on the three summands with distinct eigenvalues
	$\zeta^{-1}$, $1$, and $\zeta$, respectively. This implies that we have the 
	same decomposition in $W$, even though the summands may be smaller than in $\fW$. 
	This gives us (b) and (d). For (c), note that $\om$ maps $\im\chi$ onto $\im\om$ 
	isomorphically by (a), so $\im\om$ and $\im\chi$ have the same dimension.\qed

	\section{Solid subalgebras} \label{solid}
	
	In this section we consider an algebra $A$ of Jordan type half defined over $\F$ and 
	its subalgebra $J=\lla a,b\rra$ generated by two primitive axes $a$ and $b$
	from~$A$.
	
	\begin{defn}
		A $2$-generated subalgebra $J=\lla a,b\rra$ is \emph{solid} if every  
		idempotent $c\in J$,  $c\neq 0,1_J$, is a primitive axis of Jordan type half in $A$.
	\end{defn}
	
	Theorem \ref{main solid} from the introduction asserts that every $J$ is solid unless $\al=(a,b)\notin\{0,\frac{1}{4},1\}$. We now commence proving 
	this statement. The proof involves a number of steps. First of all, we note that $\F$ can be extended, 
	if necessary, and so, without loss of generality, we assume in the remainder of this 
	section that $\F$ is algebraically 
	closed. Throughout the proof, we also assume that $J=\lla a,b\rra$ is as above, and 
	$\al=(a,b)\notin\{0,\frac{1}{4},1\}$. In particular, by Theorem \ref{all ideals}, $J$ 
	is simple of dimension $3$, and the Frobenius form is non-degenerate on $J$. We let 
	$W=J^\perp$ be the orthogonal complement of $J$ in $A$. As before, $A=J\oplus W$ and $W$ 
	is a $J$-module.
	
	Let $\rho=\tau_a\tau_b$. We next discuss the eigenspaces of $\rho$ in $J$ and $W$.
	Note that since we assumed $\F$ algebraically closed, it contains the roots $\zeta$ 
	and $\zeta^{-1}$ of the polynomial $x^2-(4\al-2)x+1$. In particular, $\fJ(\al)$ contains 
	elements nilpotent $\fs$ and $\ft$. Set $s=\phi(\fs)$ and $t=\phi(\ft)$, where 
	$\phi:\fJ(\al)\to J$ is the unique 
	isomorphism sending $\fa$ to $a$ and $\fb$ to $b$. By Corollary \ref{on R}, $\rho$
	acts on $\la s\ra$, $\la 1_J\ra$, and $\la t\ra$ as $\zeta^{-2}$, $1$, and $\zeta^2$, 
	respectively. We note that these three numbers are distinct, except that 
	we may have $\zeta^{-2}=\zeta^2=-1$ when $\al=\frac{1}{2}$.
	
	Now we turn to $W$. Let $\om=\ad_s|_W$ and $\chi=\ad_t|_W$. We have the following statement 
	generalising Proposition \ref{3-gen W} to our arbitrary algebra $A$ 
	of Jordan type half.
	
	\begin{prop} \label{W}
		We have the following:
		\begin{enumerate}
			\item[$(a)$] $\chi\om$ acts as $\frac{1}{16}$ on $\im\om$ and, symmetrically, $\om\chi$ 
			acts as $\frac{1}{16}$ on $\im\chi$;
			\item[$(b)$] $W=\im\om\oplus(\ker\om\cap\ker\chi)\oplus\im\chi$;
			\item[$(c)$] $\rho$ acts as $\zeta^{-1}$ on $\im\om$, as $1$ on $\ker\om\cap\ker\chi$, 
			and as $\zeta$ on $\im\chi$.
		\end{enumerate}
	\end{prop}
	
	\pf Let $B$ be a basis of axes in $A$, which exists by Theorem \ref{basis of axes}. Then 
	$W$ is spanned by $B^\circ=\{c^\circ\mid c\in B\}$, where $c^\circ$ 
	is the projection of $c$ onto $W$.
	
	Focussing on a particular axis $c\in B$ and applying Proposition \ref{3-gen W} 
	to the subalgebra $\lla a,b,c\rra$, we obtain that $c^\circ=c_s+c_d+c_t$, where 
	$c_s\in\im\om$ is a $\zeta^{-1}$-eigenvector of $\rho$, $c_d\in\ker\om\cap\ker\chi$ is a 
	$1$-eigenvector of $\rho$, and $c_t\in\im\chi$ is a $\zeta$-eigenvector of $\rho$. Since 
	the elements $c^\circ$, $c\in B$, span $W$, we deduce that $W=W_s\oplus D\oplus W_t$, 
	where $W_s=\la c_s\mid c\in B\ra$, $D=\la c_d\mid c\in B\ra$ and $W_t=\la c_t\mid c\in B\ra$, 
	and the three summands, $W_s$, $D$, and $W_t$, are the $\zeta^{-1}$-, $1$-, and 
	$\zeta$-eigenspaces of $\rho$, respectively. 
	
	Furthermore, $W_s\subseteq\im\om$, $D\subseteq\ker\om\cap\ker\chi$, and $W_t\subseteq\im\chi$. 
	Again, looking at $\lla a,b,c\rra$, we observe that $(c^\circ)^\om$ is a $\zeta^{-1}$-eigenvector 
	of $\rho$, i.e., $(c^\circ)^\om\in W_s$. This means that $\im\om=\la(c^\circ)^\om\mid c\in 
	B\ra\subseteq W_s$. Hence $\im\om=W_s$ and, symmetrically, $W_t=\im\chi$. 
	
	Again, looking at $\lla a,b,c\rra$, $\chi\om$ sends each $c_s$ to $\frac{1}{16}c_s$. Hence 
	$\chi\om$ acts as $\frac{1}{16}$ on $W_s$. Symmetrically, $\om\chi$ acts as $\frac{1}{16}$ 
	on $W_t$. This gives (a) since $W_s=\im\om$ and $W_t=\im\chi$.
	
	For (b), note that $\om^2$ kills every $c^\circ$, $c\in B$, and so $\om^2=0$ and, symmetrically, 
	$\chi^2=0$. Hence $W_s\subseteq\ker\om$ and $W_t\subseteq\ker\chi$. Hence $W_s\oplus 
	D\subseteq\ker\om$ and $D\oplus W_t\subseteq\ker\chi$. As $\om\chi$ acts as $\frac{1}{16}$ 
	on $W_s$, we have that $\ker\om\cap W_t=0$. Taking into account that $W=W_s\oplus D\oplus W_t$,
	we conclude that $\ker\om=W_s\oplus D$ and, symmetrically, $\ker\chi=D\oplus W_t$. In particular, 
	$D=\ker\om\cap\ker\chi$. This completes the proof of (b) and (c).\qed
	
	\medskip
	We note that, since $\al\neq 0,1$, we have that $1$, $\zeta^{-1}$, and 
	$\zeta$ are pairwise distinct. Furthermore, $\zeta^{-2}$ is distinct from $\zeta^{-1}$ and 
	$1$. Since $\al\neq\frac{1}{4}$, we also have that $\zeta^{-2}\neq\zeta$. (Indeed, 
	$\zeta^{-2}=\zeta$ means that $\zeta$ is a cubic root of unity, i.e. a root of $x^2+x+1$. 
	Since $\zeta$ is, by definition, a root of $x^2-(4\al-2)x+1$, we see that $\zeta$ is a cubic 
	root exactly when $\al=\frac{1}{4}$.) Symmetrically, $\zeta^2$ is distinct from 
	$\zeta^{-1}$, $1$ and $\zeta$. Note that we may have that $\zeta^{-2}=\zeta^2$, which happens 
	exactly when $\al=\frac{1}{2}$. Hence we can claim that $W_s=\im\om$ is the 
	$\zeta^{-1}$-eigenspace of $\rho$ in the entire algebra $A$. Symmetrically, $W_t=\im\chi$ is 
	the $\zeta$-eigenspace of $\rho$ in $A$. Also, $\la 1_J\ra\oplus D$ is the $1$-eigenspace of 
	$\rho$ in $A$, but we cannot in general claim that $\la s\ra$ and $\la t\ra$ are 
	the $\zeta^{-2}$- and $\zeta$-eigenspaces of $\rho$ in $A$.
	
	We now set $A_0:=\la 1_J\ra\oplus D$, $A_{-1}:=W_s$, $A_1:=W_t$, $A_{-2}:=\la s\ra$, 
	$A_2:=\la t\ra$, and $A_k=0$ for all other $k\in\Z$. Clearly, we have that 
	$$A=A_{-2}\oplus A_{-1}\oplus A_0\oplus A_1\oplus A_2.$$
	
	We now claim that this decomposition is in fact a $\Z$-grading of $A$. 
	
	\begin{prop} \label{Z grading}
		We have that $A_nA_m\subseteq A_{n+m}$ for all $n,m\in\Z$.
	\end{prop}
	
	\pf Since $A_n=0$ or $A_m=0$ unless $n,m\in\{-2,-1,0,1,2\}$, we have only a finite number 
	of pairs $(n,m)$ to consider. So let $u\in A_n$ and $v\in A_m$ with 
	$n,m\in\{-2,-1,0,1,2\}$. Furthermore, by symmetry, it suffices to consider the following 
	pairs $(n,m)$: $(-2,-2)$, $(-2,-1)$, $(-2,0)$, $(-2,1)$, $(-2,2)$, $(-1,-1)$, $(-1,0)$, 
	$(-1,1)$, and $(0,0)$.
	
	If $n+m\in\{-1,0,1\}$ then we can use the eigenspace argument. We note that $u^\rho=
	\zeta^n u$ and $v^\rho=\zeta^m v$. Hence $(uv)^\rho=u^\rho v^\rho=\zeta^m u\zeta^n v=
	\zeta^{n+m} uv$, and so $uv$ must be in the $\zeta^{n+m}$-eigenspace, which coincides 
	with $A_{n+m}$. So the claim is true in the cases $(n,m)=(-2,1)$, $(-2,2)$, $(-1,0)$, 
	$(-1,1)$, and $(0,0)$.
	
	If $(n,m)=(-2,-2)$ then $u$ and $v$ are multiples of $s$, and so $uv=0$ since $s^2=0$. 
	This shows that $uv\in A_{m+n}=A_{-4}=0$, and so the claim holds in this case.
	
	If $(n,m)=(-2,-1)$ then $u$ is a multiple of $s$ and $v\in W_s$. This implies that 
	$uv\in\om(W_s)=0$, since $W_s\subseteq\ker\om$. Thus, again $uv=0$ and it 
	is contained in $A_{m+n}=A_{-3}=0$.
	
	If $(n,m)=(-2,0)$ then $u$ is a multiple of $s$ and $v\in\la 1_J\ra\oplus D$. 
	So $v=\pi 1_J+d$ for some $\pi\in\F$ and $d\in D$. Then $uv=\pi u1_J+ud=\pi u+0$, 
	since $d\in\ker\om$. Thus, $uv=\pi u\in\la s\ra=A_{-2}=A_{n+m}$.
	
	Finally, if $(n,m)=(-1,-1)$ then $uv$ is contained in the $\zeta^{-2}$-eigenspace. 
	If $\al\neq\frac{1}{2}$ then this eigenspace coincides with $\la s\ra=A_{-2}$, 
	and hence the claim holds. If $\al=\frac{1}{2}$ then $\zeta^{-2}=-1=\zeta^2$, and so 
	the $\zeta^{-2}$-eigenspace coincides with $\la s\ra\oplus\la t\ra$. So 
	$uv=\pi s+\epsilon t$ for some $\pi,\epsilon\in\F$. On the other hand, 
	$0=(u,0)=(u,vs)=(uv,s)=(\pi s+\epsilon t,s)=\pi(s,s)+\epsilon(s,t)=\pi 0+
	\epsilon\frac{1}{4}=\frac{\epsilon}{4}$. Thus, $\epsilon=0$ and hence $uv\in 
	\la s\ra=A_{-2}=A_{n+m}$. This completes the last case, and so the claim is 
	proven.\qed
	
	\medskip
	This immediately implies the following.
	
	\begin{cor}
		For $\xi\in\F$, $\xi\neq 0$, the linear map $\phi_\xi:A\to A$ that acts as $\xi^n$ 
		on each $A_n$, $n\in\{-2,-1,0,1,2\}$, is an automorphism of $A$.
	\end{cor}
	
	Let $\Phi=\{\phi_\xi\mid\xi\in\F,\xi\neq 0\}$. This is a subgroup of $\Aut(A)$.
	
	\begin{prop}
		The group $\Phi$ acts transitively on the set of idempotents from $J$, other than $0$ 
		and $1_J$.
	\end{prop}
	
	\pf Consider idempotents $u=\dl s+\frac{1}{2}1_J+\dl^{-1}t$ and $v=\epsilon s+
	\frac{1}{2}1_J+\epsilon^{-1}t$ with $\dl,\epsilon\in\F\setminus\{0\}$. Let 
	$\xi\in\F$ such that $\xi^2=\frac{\dl}{\epsilon}$. Then $u^{\phi_\xi}=
	\dl s^{\phi_\xi}+\frac{1}{2}1_J^{\phi_\xi}+\dl^{-1}t^{\phi_\xi}=\dl\xi^{-2}s+
	\frac{1}{2}\xi^01_J+\dl^{-1}\xi^2t=\epsilon s+\frac{1}{2}1_J+\epsilon^{-1}t=v$. 
	So the claim holds.\qed
	
	\medskip
	Since the axes $a$ and $b$ are among the idempotents in this proposition, all 
	these idempotents are axes of $A$, and so $B$ is solid, which completes the proof of 
	Theorem \ref{main solid}.

	\section{Characteristic zero} \label{char zero}
	
	When the field $\F$ is of characteristic zero, we can remove some of the exceptional 
	configurations using algebraic geometry arguments. We first make some general observations.
	
	Let $A$ be a finite-dimensional algebra over of field $\F$. Say, $n=\dim(A)$.
	
	\begin{prop}
		Idempotents in $A$ form an affine algebraic set.
	\end{prop}
	
	\pf The condition $u^2=u$ amount to $n$ quadratic equations in the coordinates of $u\in A$ with respect to an arbitrary basis $B$ of $A$.\qed
	
	\medskip
	More interestingly, we also have the following.
	
	\begin{prop}
		If $\cF\subseteq\F$ is a (finite) fusion law then the set of $\cF$-axes in $A$ is also an affine algebraic set, a subset of the algebraic set of idempotents.  
	\end{prop}
	
	\pf We need to express the conditions on primitive $\cF$-axes in terms of polynomials in the coordinates of the axis. Let $f(x)\in\F[x]$ be defined as $f(x)=\prod_{\lm\in\cF}(x-\lm)$. 
	Furthermore, for each subset $\cH\subseteq\cF$, we set $f_\cH(x)=\prod_{\lm\in\cH}(x-\lm)$, so that 
	$f(x)=f_\cF(x)$. We also let $g_\lm(x)=\frac{f(x)}{x-\lm}=f_{\cF\setminus\{\lm\}}(x)$ for $\lm\in\cF$.
	
	Then the condition that $\ad_a$, for an idempotent $a$, is semi-simple with the entire spectrum within $\cF$ can be expressed by the equality $f(\ad_a)=0$, and this clearly is equivalent to $n^2$ (since we can view $\ad_a$ as an $n\times n$ matrix) polynomial conditions in the $n$ coordinates of $a$ with respect to the basis $B$.
	
	For the fusion law, note that the $\lm$-eigenspace of $\ad_a$ (that already satisfies the above conditions) can be expressed as follows: $A_\lm(a)=g_\lm(\ad_a)A$. Indeed, $g_\lm(\ad_a)$ acts as zero on all eigenspaces $A_\mu(a)$ with $\mu\neq\lm$, and it acts on $A_\lm(a)$ as the scalar $g_\lm(\lm)=\prod_{\mu\in\cF\setminus\{\lm\}}(\lm-\mu)\neq 0$. Furthermore, the set $\{g_\lm(\ad_a)v\mid b\in B\}$ is a spanning set for $A_\lm(a)$, since $B$ is a basis of $A$.
	
	Hence, the fusion law conditions can be expressed as follows: for all $\lm,\mu\in\cF$ and all $b,b'\in B$, we must have that $f_{\lm\ast\mu}(\ad_a)((g_\lm(\ad_a)b)(g_\mu(\ad_a)b'))=0$. (This is because $f_{\lm\ast\mu}(\ad_a)$ acts as zero on $A_{\lm\ast\mu}(\ad_a)$ and nowhere else.) All these conditions end up being a collection of polynomial equations in the $n$ coordinates of $a$ with respect to $B$.\qed
	
	\medskip
	We note that this argument is based on the ideas of the universal algebra construction from \cite{hrs1}. Here is an additional observation from there.
	
	\begin{prop} \label{axis variety}
		If $A$ admits a Frobenius form then the set of all primitive axes satisfying $(a,a)\neq 0$ is 
		also an algebraic set, a subset of the algebraic set of all $\cF$-axes.
	\end{prop}
	
	\pf Indeed, the additional condition that $a$ is a primitive axis can be expressed as $g_1(\ad_a)(b-\frac{(b,a)}{(a,a)}a)=0$, and clearly, this needs to be accompanied by the condition $(a,a)z=1$, where $z$ is an additional indeterminate. All these conditions are polynomial.\qed
	
	\medskip
	When $A$ is of Jordan type half, it admits a Frobenius form such that all primitive axes have length $(a,a)=1$, and so the extra variable $z$ is not needed in this case.
	
	Let us now apply these general observations for the purposes of this paper. We want to show that at least in characteristic zero the additional $2$-generated algebras, arising for $\al=1$ and $\al=0$ are solid, that is, all the primitive idempotents in them are primitive axes of $A$.
	
	We first consider the case $\al=1$. The algebras $\fJ(1)=\lla\fa,\fb\rra$ and its factor $\overline{\fJ(1)}$ are baric, that is, they have radical of codimension $1$. These algebras do not have an identity element.
	
	\begin{lem}
		The affine algebraic set of non-zero idempotents in $\fJ(1)$ and $\overline{\fJ(1)}$ is of dimension $1$. 
	\end{lem} 
	
	\pf In the case of $A=\fJ(1)$, we set $u=\lm_1\fa+\lm_2\fb+\lm_3\sg$ and deduce:
	$u^2=(\lm_1\fa+\lm_2\fb+\lm_3\sg)(\lm_1\fa+\lm_2\fb+\lm_3\sg)=\lm_1^2\fa+2\lm_1\lm_2(\frac{1}{2}\fa+\frac{1}{2}\fb+\sg)+\lm_2^2\fb$ since $\sg$ annihilates the entire algebra when $\al=1$. So $u^2=(\lm_1^2+\lm_1\lm_2)\fa+(\lm_2^2+\lm_1\lm_2)\fb+2\lm_1\lm_2\sg$. Thus, the equality $u^2=u$ gives us the equations:
	$$\lm_1(\lm_1+\lm_2-1)=0,$$
	$$\lm_2(\lm_1+\lm_2-1)=0,$$
	and
	$$2\lm_1\lm_2=\lm_3.$$
	
	From these we see that either $\lm_1=\lm_2=\lm_3=0$, or $\lm_1+\lm_2=1$ and $\lm_3=2\lm_1\lm_2$. We conclude that the algebraic set of idempotents has two irreducible components: $\{0\}$ and the $1$-parameter components of non-zero idempotents.
	
	The calculations for $\overline{\fJ(1)}$ are even simpler, since there is no $\sg$ and $\lm_3$, and we skip them.\qed
	
	\begin{prop} \label{exceptional 1}
		Suppose $A$ is an algebra of Jordan type half over a field $\F$ is of characteristic zero. If $a,b\in A$ are primitive axes with $\al=(a,b)=1$ then all non-zero idempotents in $J=\lla a,b\rra$ are primitive axes in $A$. That is, $J$ is solid.
	\end{prop}
	
	\pf We have seen above that the conditions assuring that $u$ is a primitive axis are polynomial, and so the set of those idempotents $u\in J$ that are axes in the entire $A$ form a sub algebraic set of the algebraic set of non-zero idempotents of $J$. 
	
	Since the latter is irreducible of dimension $1$, its every proper sub algebraic set has dimension zero and hence is finite. However, $\rho=\tau_a\tau_b$ in this case has infinite order and so infinitely many idempotents $a^{\rho^i}$ and $b^{\rho^i}$ are  primitive axes. Hence this sub algebraic set has dimension $1$ and it coincides with the set of all non-zero idempotents from $J$.\qed
	
	\medskip
	Next we deal with the case $\al=0$. In this case, the $2$-generated algebra $J=\lla a,b\rra$ with $(a,b)=0$ is either the $3$-dimensional algebra $\fJ(0)$ or its factor $\bar\fJ(0)\cong 2B$. The latter is trivially solid, as its only primitive idempotents are $a$ and $b$. So we focus on $J\cong\fJ(0)$.
	
	\begin{prop} \label{value 0}
		The affine algebraic set of primitive non-zero idempotents in $\fJ(0)$ consists of two disjoint irreducible components of dimension one. Furthermore, the two generators, $\fa$ and $\fb$, of $\fJ(0)$ belong to different components.
	\end{prop}
	
	\pf Let $u=\lm_1\fa+\lm_2\fb+\lm_3\sg$. Then we deduce:
	$u^2=(\lm_1\fa+\lm_2\fb+\lm_3\sg)(\lm_1\fa+\lm_2\fb+\lm_3\sg)=\lm_1^2\fa+2\lm_1\lm_2(\frac{1}{2}\fa+\frac{1}{2}\fb+\sg)+\lm_2^2\fb+2\lm_1\lm_3(-\frac{1}{2}\fa)+2\lm_2\lm_3(-\frac{1}{2}\fb)+\lm_3^2(-\frac{1}{2}\sg)=(\lm_1^2+\lm_1\lm_2-\lm_1\lm_3)\fa+(\lm_1\lm_2+\lm_2^2-\lm_2\lm_3)\fb+(2\lm_1\lm_2-\frac{1}{2}\lm_3^2)\sg$. 
	
	Hence we have three equations:
	$$\lm_1(\lm_1+\lm_2-\lm_3-1)=0,$$
	$$\lm_2(\lm_1+\lm_2-\lm_3-1)=0,$$
	and
	$$4\lm_1\lm_2-\lm_3^2-2\lm_3=0.$$
	
	From here, we see that either $\lm_1=\lm_2=0$ and $\lm_3=0$ or $-2$, which gives us the non-primitive zero and identity idempotents, or $\lm_3+1=\lm_1+\lm_2$. Squaring this, we get 
	$\lm_3^2+2\lm_3+1=\lm_1^2+2\lm_1\lm_2+\lm_2^2$, and so the third equation above is equivalent to 
	$4\lm_1\lm_2-\lm_1^2-2\lm_1\lm_2-\lm_2^2+1=0$, which simplifies to $1-(\lm_1-\lm_2)^2=0$. This factorises as $(1-\lm_1+\lm_2)(1+\lm_1-\lm_2)=0$. Thus we have two components on primitive idempotents, one given by the linear equations 
	$$\left\{
	\begin{array}{rrrrrrr}
		\lm_1&-&\lm_2&&&=&1\\
		\lm_1&+&\lm_2&-&\lm_3&=&1
	\end{array}
	\right.$$
	and the other by the linear equations
	$$\left\{
	\begin{array}{rrrrrrr}
		\lm_1&-&\lm_2&&&=&-1\\
		\lm_1&+&\lm_2&-&\lm_3&=&1
	\end{array}
	\right.$$
	
	Clearly, this means we have two disjoint $1$-dimensional irreducible components. Furthermore, manifestly, $\fa$ (for which $(\lm_1,\lm_2,\lm_3)=(1,0,0)$) is in the first component and $\fb$ ($(\lm_1,\lm_2,\lm_3)=(0,1,0)$) is in the second component.\qed
	
	\begin{prop} \label{exceptional 0}
		Suppose $A$ is an algebra of Jordan type half over a field $\F$ is of characteristic zero. If $a,b\in A$ are primitive axes with $\al=(a,b)=0$ then all non-zero, non-identity idempotents in $J=\lla a,b\rra$ are primitive axes in $A$. That is, $J$ is solid.
	\end{prop}
	
	\pf If $J\cong 2B$ then the claim is true, since $a$ and $b$ are the only primitive idempotents in $J$.
	
	Now let us assume that $J$ is of dimension $3$. As above, by Proposition \ref{axis variety}, primitive axes of $A$ contained in $J$ form a sub algebraic set $Y$ of the affine algebraic set $Z$ of all primitive idempotents of $J$. According to Proposition \ref{value 0}, $Z$ consists of two $1$-dimensional irreducible components. We claim that both components are contained in $Y$, that is, $X=Y$. Indeed, by the same Proposition \ref{value 0}, $a$ and $b$ belong in different components of $Z$. Hence each component contains $x=a$ or $b$ and also infinitely many other axes $x^{\la\tau_a,\tau_b\ra}$. Hence the intersection of the component with $Y$ cannot be of dimension $0$, and so the entire component must be fully contained in $Y$.\qed
	
	\medskip
	Propositions \ref{exceptional 1} and \ref{exceptional 0} establish Theorem \ref{exceptional}. We now have that, in characteristic $0$, essentially all $2$-generated subalgebras in $A$ are solid. Namely, we can state the following.
	
	\begin{cor} 
		Let $\F$ be of characteristic $0$, $A$ be an algebra of Jordan type half and $J=\lla a,b\rra$ for primitive axes $a,b\in A$. If $J$ is not solid then $(a,b)=\frac{1}{4}$ and $J\cong 3C(\frac{1}{2})$ contains exactly three primitive axes from $A$, namely $a$, $b$, and $c=a^{\tau_b}=b^{\tau_a}$.
	\end{cor}
	
	Among the solid lines, if $\F$ is infinite, only the algebra $J=\lla a,b\rra\cong 2B$, where $ab=0$ and $a$ and $b$ are the only primitive axes, contains finitely many axes. 
	This allows us to determine all algebras $A$ with finitely many primitive axes (equivalently, a finite automorphism group). In other words, we now establish Theorem \ref{main corollary}.
	
	In the remainder of this section, $A$ is an algebra of Jordan type half, with a closed (finite) generating set of axes $X$, such that for all $a,b\in X$, either $(a,b)=0$ and $ab=0$ (and so, $\lla a,b\rra\cong 2B$) or $(a,b)=\frac{1}{4}$ (and so $\lla a,b\rra\cong 3C(\frac{1}{2})$). In the latter case, we have that $ab=\frac{1}{4}(a+b-c)$, where 
	$c=a^{\tau_b}=b^{\tau_a}$.
	
	Without loss of generality, we may assume that for each $a\in X$ there is at least one $b\in X$ such that $ab\neq 0$. (Otherwise, $\la a\ra$ is a direct summand of $A$ 
	by Seress Lemma.) Under this additional assumption, $\tau_a\neq 1$ for all $a\in X$. 
	
	\begin{prop}
		The Miyamoto group $G=\Miy(A)=\la\tau_a\mid a\in X\ra$ is a $3$-transposition group with respect to the normal set $D=\{\tau_a\mid a\in X\}$.
	\end{prop}
	
	\pf Clearly, when $ab=0$, we have that the corresponding involutions, $\tau_a$ and $\tau_b$, commute since $\tau_a^{\tau_b}=\tau_{a^{\tau_b}}=\tau_a$. So, in this case, 
	$(\tau_a\tau_b)^2=[\tau_a,\tau_b]=1$. Thus, we just need to consider the case where $(a,b)=\frac{1}{4}$. In this case, set $H=\la\tau_a,\tau_b\ra$ and $z=(\tau_a\tau_b)^3$. 
	Clearly, $H$ acts on the set $\{a,b,c\}$, where $c=a^{\tau_b}$, and furthermore, $z$ acts on this set trivially, i.e., it fixes all three axes. In particular, 
	$(\tau_a)^z=\tau_{a^z}=\tau_a$ and, similarly, $\tau_b^z=\tau_b$. That is, $z$ is in the centre of $H$. However, $H$ is a dihedral group and $\tau_a$, say, inverts 
	all powers of $\rho=\tau_a\tau_b$. This implies that either $z=1$ (which is the conclusion we need) or $z$ has order $2$. However, in the latter case, the order of 
	$\rho$ is $6$ and $H\cong D_{12}$. This leads to a contradiction, since in $D_{12}$ the two generating involutions, $\tau_a$ and $\tau_b$, are not conjugate, while
	in our situation $\tau_b=\tau_a^{\tau_c}$ (and $\tau_c=\tau_{b}^{\tau_a}\in H$). Thus, indeed, $z=1$ and $|\tau_a\tau_b|=|\rho|=3$.\qed
	
	\medskip
	This allows us to establish the claim in Theorem \ref{main corollary}. Note that, by Theorem B from \cite{hss1}, if $\tau_a=\tau_b$ for distinct axes 
	$a,b\in X$ then the connected component $B$ of $A$ containing $a$ and $b$ is a spin factor Jordan algebra. In particular, $a+b=1_B$. However, if $c\in X$ is such 
	that $(a,c)=\frac{1}{4}$ we then have that $(b,c)=(1_B-a,c)=(1_B,c)-(a,c)=1-\frac{1}{4}=\frac{3}{4}$, which is a contradiction. Hence, the map $a\mapsto\tau_a$ 
	is a bijection between $X$ and $D$.
	
	\begin{prop}
		There exists a surjective homomorphism from $M_{\frac{1}{2}}(G,D)$ onto $A$; i.e., $A$ is isomorphic to a Matsuo algebra or a factor of Matsuo algebra.
	\end{prop}
	
	\pf Let $M=M_{\frac{1}{2}}(G,D)$ and $\phi:M\to A$ be the linear map sending each $\tau_a$ back to $a$. (Recall that $D=\{\tau_a\mid a\in X\}$ is a basis of $M$.) To check that $\phi$ is an algebra homomorphism, it suffices to show that
	$$\phi(\tau_a\cdot\tau_b)=\phi(\tau_a)\phi(\tau_b)=ab,$$ 
	for all $a,b\in X$. And indeed, if $ab=0$ then, by the definition of Matsuo algebra, $\tau_a\cdot\tau_b=0$ and so both sides above are zero. In the second case, 
	where $(a,b)=\frac{1}{4}$, we also have that $\phi(\tau_a\cdot\tau_b)=\phi(\frac{1}{4}(\tau_a+\tau_b-\tau_a^{\tau_b}))=\frac{1}{4}(a+b-c)=ab$, where $c=a^{\tau_b}$
	and so $\tau_a^{\tau_b}=\tau_c$.\qed
	
	\medskip
	To conclude this section, we note that, in the positive characteristic, the algebraic geometry method we used in this section, falls apart: indeed, in this case, for both $(a,b)=1$ and $(a,b)=0$, we have that $a^{\la\tau_a,\tau_b\ra}$ and $b^{\la\tau_a,\tau_b\ra}$ are finite of size $p=\mathrm{char}(\F)$ and we cannot claim that the set of primitive axes is of positive dimension.

	\section{Example} \label{example}
	
	De Medts and Rehren \cite{dmr} determined Matsuo algebras with $\eta=\frac{1}{2}$ that are Jordan algebras. (There is also a recent article \cite{gms} which identifies in zero characteristic all Matsuo algebras having non-zero Jordan factor algebras.) In such an algebra $M$, all idempotents, even the non-primitive ones, satisfy the Jordan law of type half. Hence in $M$ all subalgebras are trivially solid. Below we provide an example found by Gorshkov and Staroletov within the project on the classification of $4$-generated algebras of Jordan type half by their diagrams. In this example, the Matsuo algebra is not Jordan, but it nevertheless contains solid lines.
	
	First we need the following lemma that provides a useful way of checking whether a given subalgebra in a Matsuo algebra is solid.
	
	\begin{lem} \label{a-prime}
		Let $M=M_{\frac{1}{2}}(G,D)$ be the Matsuo algebra of the $3$-transposition group $(G,D)$. Let $a$ and $b$ be collinear points in the Fischer space of $(G,D)$. Then the subalgebra $J=\lla a,b\rra$ of $M$ is solid if and only if $a'=1_J-a$ is an axis in $M$.
	\end{lem}
	
	\pf We note that, by Lemma 9 in \cite{gg}, $a'$ is a primitive idempotent with spectrum $\{1,0,\frac{1}{2}\}$. So the only issue is whether the Jordan fusion law is satisfied for $a'$. 
	
	First of all, if $J$ is solid, then every primitive idempotent in $J$ is an axis in $M$, and 
	in particular, $a'$ is an axis in $M$.
	
	Conversely, suppose that $a'$ is an axis in $M$. Let $J'=\lla a',b\rra$. Clearly, $J'\subseteq B$. Furthermore, $(a',b)=(1_J-a,b)=(1_J,b)-(a,b)=(1_J,b^2)-\frac{1}{4}=(1_Jb,b)-\frac{1}{4}=
	(b,b)-\frac{1}{4}=1-\frac{1}{4}=\frac{3}{4}$. According to Theorem \ref{main solid}, $J'$ is solid. Furthermore, $J'$ is simple of dimension $3$, and so $J'=J$. Thus, $J$ is solid.\qed
	
	\begin{exa} \label{vertical}
		Let $G=3^3:S_4$ be the $3$-transposition group denoted in \cite{hs} as $3^{\bullet 1}:S_4$. Here $D$ is the class of involutions that generates $G$. Then $|D|=18$. Let $\pi$ be the natural homomorphism from $G$ onto $\bar G\cong S_4$. Then $\bar D=\pi(D)$ is of size $6$. The Fischer space on $D$ has two types of lines: (a) lines $L$, such that $\bar L=\pi(L)$ is of size $3$, that is, $\bar L$ is in a bijection with $L$; and (b) lines $L$, such that $\bar L$ is of size $1$. The latter type of lines are the fibres of the map $\pi$ restricted to $D$. We will call the lines in (a) \emph{horizontal} and the lines in (b) \emph{vertical}. These two sets of lines are the orbits of $G$ on the set of all lines.
		
		It was verified computationally, using Lemma \ref{a-prime}, that in this example vertical lines are solid, while the horizontal lines are not solid. In particular, $M$ cannot be Jordan, since it contains non-solid lines. We note that this algebra $M$ has a $2$-dimensional ideal $I$, the radical of $M$, and $M/I$ is a Jordan algebra of dimension $16$.
	\end{exa}
	
	This example leads to a number of questions\footnote{Jari Desmet informed the authors that he obtained a positive answer to the first two questions}.
	
	\begin{que}
		Is it true that the vertical lines are solid for all Matsuo algebras $M_{\frac{1}{2}}(3^{\bullet 1}S_n)$?
	\end{que}
	
	Note that $S_n=W(A_{n-1})$ is the Weyl group of type $A_{n-1}$. Hence the above question can be generalised further.
	
	\begin{que}
		Is it true that vertical lines are solid in all Matsuo algebras $M_{\frac{1}{2}}(3^{\bullet 1}W(X_n))$ for all simply laced diagrams $X=A$, $D$, and $E$? 
	\end{que}
	
	The ultimate question is:
	
	\begin{que}
		Is it possible to determine all Matsuo algebras and their factors that admit solid lines. 
	\end{que}
	
	It would also be interesting to generalise the above construction to groups $3^{\bullet 1}W(X_n)$ for non simply laced diagrams. Here, of course, we cannot talk about Matsuo algebras, and so this would also require defining suitable axial algebras for these groups.

	\Addresses

\begin{thebibliography}{WWW}
		
		\bibitem[C]{c}
		J.H.~Conway,  A simple construction for the Fischer-Griess monster group, {\it Ivent. Math.}, {\bf 79} (1985), 513--540.
		
		\bibitem[DMPSVC]{dmpsvc}
		T.~De Medts, S.~Peacock, S.~Shpectorov, M.~Van Couwenberghe,
		Decomposition algebras and axial algebras, {\it J. Algebra}, {\bf 556} (2020), 287--314.
		
		\bibitem[DMR]{dmr}
		T.~De Medts, F.~Rehren, Jordan algebras and $3$-transposition groups, {\it J. Algebra}, {\bf 478} (2017), 318--340.
		
		\bibitem[DMRS]{dmrs}
		T.~De Medts, L.~Rowen, Y.~Segev, Primitive $4$-generated axial algebras of Jordan type,  {\it Proc. Amer. Math. Soc.}, DOI:10.1090/proc/16571.
		
		\bibitem[GAP]{GAP}
		The GAP Group, GAP -- Groups, Algorithms, and Programming, Version 4.12.2, 2022. (https://www.gap-
		system.org)
		
		\bibitem[GG]{gg}
		I.~Gorshkov, V.~Gubarev, Quasi-definite axial algebras of Jordan type half, {\it Siberian Electronic Mathematical Reports}, {\bf 20}:2 (2023), 833--846.
		
		\bibitem[GMS]{gms}
		I.~Gorshkov, A.~Mamontov, and A.~Staroletov, On Jordan algebras that are factors of Matsuo algebras,
		manuscript, arXiv:2305.10958.
		
		\bibitem[GS]{gs} I.~Gorshkov and A.~Staroletov, On primitive $3$-generated axial 
		algebras of Jordan type, {\it J. Algebra} {\bf 563} (2020), 74–-99.
		
		\bibitem[G]{g}
		R.L.~Griess Jr., The friendly giant, {\it Ivent. Math.}, {\bf 69} (1982), 1--102.
		
		\bibitem[HRS]{hrs1} J.I.~Hall, F.~Rehren, and S.~Shpectorov, Universal Axial Algebras and 
		a Theorem of Sakuma, {\it J. Algebra}, {\bf 421} (2015), 394--424.
		
		\bibitem[HRS2]{hrs2} J.I.~Hall, F.~Rehren, and S.~Shpectorov, Primitive axial 
		algebras of Jordan type, {\it J. Algebra} {\bf 437} (2015), 379-–115.
		
		\bibitem[HSS]{hss1} J.I.~Hall, Y.~Segev, S.~Shpectorov, Miyamoto involutions in axial algebras of Jordan type half, {\it Israel J. Math.} {\bf 223} (2018), 261--308.
		
		\bibitem[HSS2]{hss2} J.I.~Hall, Y.~Segev, and S.~Shpectorov, 
		On primitive axial algebras of Jordan type, {\it Bull. Inst. Math. Acad. Sin. (N.S.)} {\bf 13}:4 (2018), 397--409.
		
		\bibitem[HS]{hs} J.I.~Hall, S.~Shpectorov, The spectra of finite $3$-transposition groups, {\it Arab. J. Math.} {\bf 10}:3 (2021), 611--638. 
		
		\bibitem[I]{i} A.A.~Ivanov, The Monster Group and Majorana Involutions,
		{\it Cambridge Tracts in Mathematics} {\bf 176}, Cambridge University Press,
		2009, 266 pp.
		
		
		\bibitem[KMS]{kms} 
		S.M.S.~Khasraw, J.~McInroy, S.~Shpectorov, 
		On the structure of axial algebras, {\it Trans. Amer. Math. Soc.} {\bf373}:3 (2020), 2135--2156.
		
		\bibitem[M]{m} A.~Matsuo, 3-transposition groups of symplectic type and 
		vertex operator algebras (version 1), manuscript, November 2003 (available as
		arXiv.math/0311400v1).
		
		\bibitem[MS]{ms} J.~McInroy, S.~Shpectorov, Axial algebras of Jordan and Monster type, arXiv:2209.08043.
		
	\end{thebibliography}
\end{document}